%% file: sphere.tex
\documentclass[12pt]{article}
\usepackage{a4,amsthm,amsmath,amssymb}
\usepackage{graphicx,color}

\newcommand{\bx}{\hspace*{\fill}$\square$}

\begin{document}
\noindent
{\Large \bf On the set of complex points of a 2-sphere}\\
\vskip0.4cm
\noindent
{\large \bf Nikolay Shcherbina}\\

\noindent
{\small Department of Mathematics, University of Wuppertal, 42119  Wuppertal, Germany}\\
({\small e-mail: \tt{shcherbina@math.uni-wuppertal.de}})\\

\vskip0.6cm

\noindent
{\large \bf {1. Introduction}}\\

\noindent
Let $M$ be a 2-dimensional $C^1$-smooth manifold in ${\mathbb C}^2$. 
A point $p$ on $M$ is called a {\it complex point}, if the tangent 
plane $T_pM$ to $M$ at $p$ is a complex line. Denote by ${\mathcal E}$ 
the set of all complex points on $M$. If $M$ is smooth enough and in a 
general position, then the set ${\mathcal E}$ consists of isolated points.
In this case the topology of $M$ can be described in terms of the local 
behaviour of $M$ near the points of ${\mathcal E}$ (see [L]).
The structure of the set $M$ near the points in ${\mathcal E}$ plays 
a key role in different questions of complex analysis
(see, for example, [Bi], [BK], [K], [N], [Wi] and [J]). 
In this paper we study the structure of the set ${\mathcal E}$ in the case, 
when $M$ is a 2-dimensional sphere, denoted by $S$ in what follows, embedded  
into the boundary $\partial G$ of a $C^\infty$-smooth strictly pseudoconvex  
domain $G$ in ${\mathbb C}^2$ (this case is important for applications, as  
was shown in [El] and [Er]). It is easy to see that, due to strict pseudoconvexity 
of $G$, the set ${\mathcal E}$ has no interior points in $S$. Our goal here is 
to give a more explicit description of the set ${\mathcal E}$ depending on the 
smoothness of $S$.  Recall, that a manifold is said to be of class $C^{2-}$ if 
it can be represented locally as the graph of a function that belongs to the class 
${\rm Lip}^{1,\alpha}$  for  each positive $\alpha < 1$. Our main result can now 
be formulated as follows.

\vskip0.5cm
\noindent
{\bf Theorem.} {\it Let $G$ be a strictly pseudoconvex domain in ${\mathbb C}^2$ with 
$C^\infty$-smooth boundary $\partial G$. Let $S$ be a 2-dimensional sphere embedded 
into $\partial G$. Then, depending on the smoothness of $S$, the following holds:

\vskip0.5cm
1) 
If $S$ is of class $C^2$, then the set ${\mathcal E}$ of complex points of $S$ is 
contained in a $C^1$-smooth nonclosed curve $\gamma \subset S$.

\vskip0.3cm
2)
There exists a 2-dimensional sphere $S \subset  \partial G$ of class $C^{2-}$ such 
that the set ${\mathcal E}$ contains a Jordan curve of positive 2-dimensional measure.}

\vskip0.5cm 

{\small {\it Acknowledgements.} Part of this work was done while the author was a visitor 
at the Max Planck Institute of Mathematics (Bonn) and at the Scuola Normale Superiore (Pisa). 
It is my pleasure to thank these institutions for their hospitality and excellent working 
conditions. I also express my gratitude to N. G. Kru\v{z}ilin for valuable comments.}\\

\noindent
{\large \bf {2. Proof of the first part of the theorem}}\\

\noindent
We start with an argument which goes back to Bishop [Bi] (see also [J]). Namely, 
if $p$ is a point of ${\mathcal E}$, then after a polynomial change of coordinates that 
moves $p$ to the origin we can locally represent $S$ as the disc 
$$D = \left\{(z, f(z)) \in   {\mathbb C}^2 : z \in  \triangle \right\}$$ 
with $\triangle$ being a small disc centered at the origin and $f$ being a complex  
valued $C^2$-smooth function. Moreover, in view of strict pseudoconvexity of $G$, 
after this change of coordinates the function $f$ will have the special form 
$$f(z) = \frac{1}{2} |z|^2 - \beta {\,\rm Re}{\,z^2} + o (|z|^2),\,\,\,\,\beta \geq  0$$ 
near zero. Recall, that zero is called an {\it elliptic point} if $\,0 \leq  \beta  <  \frac{1}{2}$, 
a {\it hyperbolic point} if $\,\beta  > \frac{1}{2}\,$ and a {\it parabolic point} if  
$\,\beta  = \frac{1}{2}$. Elliptic and hyperbolic points are always isolated in ${\mathcal E}$. 
In the case of a parabolic point we can use the real coordinates $z = x + iy$ and represent 
$f$ as $f(z) = y^2 + o (|z|^2)$. Hence $\,\partial_{\bar{z}}  f(z) = iy + o (|z|)\,$ and then, 
by the implicit function theorem, we obtain that
$\,\sigma  = \left\{z \in  \Delta   : {\rm Im} \,\,\partial_{\bar{z}}  f (z) = 0\right\}\, $  is a 
$C^1$-smooth curve and locally ${\mathcal E} =\left\{(z, f(z)) : z \in \sigma \,\,{\rm and} \,\, 
{\rm Re}\,\, \partial _{\bar{z}}  f(z) = 0\right\}$. Therefore, locally the set ${\mathcal E}$ 
is a closed subset of a $C^1$-smooth curve.

\vskip0.2cm 

Since the set ${\mathcal E}$ is compact, the only obstruction for ${\mathcal E}$ to be a 
subset of a nonclosed $C^1$-smooth curve $\gamma  \subset  S$ is that there is a closed 
$C^1$-smooth curve   $ \Gamma \subset  {\mathcal E}$.

\vskip0.2cm 

Assume, to get a contradiction, that such a closed $C^1$-smooth curve   $\Gamma \subset  {\mathcal E}$ exists.
Consider a complex tangential $C^\infty$-smooth closed curve
$\Gamma'$  in $\partial G$ 
(complex tangential here means that for each point $p \in \Gamma'$ the tangent line $T_p {\,\Gamma'}$ 
to $\Gamma'$ at $p$ is contained in the complex tangent plane $T_p^{\mathbb C}(\partial G) $ to 
$\partial G$ at $p$) close enough to $\Gamma$ in the $C^1$-metric. Then, using
a partition of unity along the curve $\Gamma'$, we can construct a small 
$C^1$-perturbation $S'$ of $S$ in $\partial G$ such that $\Gamma '  \subset  S'$ and each point 
in $\Gamma '$  is a complex point on $S'$. Moreover, $S'$ can be made $C^\infty$-smooth 
in a neighbourhood of $\Gamma '$. For each point $p \in \Gamma '$, consider the unit vector 
$\vec{u}(p)$ tangent to $\Gamma '$, the vector $i \vec{u} (p) \in  T_p^{\mathbb C}(\partial G) $ 
and the unit vector $\vec{n}(p) \in T_p (\partial G) $ orthogonal to $T_p^{\mathbb C}(\partial G) $ 
and such that the vectors $(\vec{u}(p), i\vec{u}(p), \vec{n}(p))$ define the positive orientation 
of $\partial G$ at the point $p$. Let $O(p)$ be the rotation of $T_p(\partial G)$ around the 
direction $\vec{u}(p)$ that transforms the vector $i \vec{u}(p)$ into the vector $\vec{n}(p)$. Using 
the tubular neighbourhood theorem (see, for example, Theorem 1.4 in [H]) we can change $S'$ in a neighbourhood
of $\Gamma '$ to get a new 2-sphere, $S'' \subset \partial G\,\,$, $C^\infty$-smooth near $\Gamma '$ such that 
$\Gamma ' \subset  S''$ and $\vec{n}(p) \in  T_p (S'')$ for each point $p \in  \Gamma '$. It is easy 
to see that $S''$ is totally real near $\Gamma '$. Then we can perturb $S''$ slightly outside a small
neighbourhood  of $\Gamma '$ to get a $C^\infty$-smooth 2-sphere $\tilde{S} \subset \partial G$ in general
position.  To finish the proof of the first part of our theorem we use an argument of Gromov (see [G], p.
342). Namely, by the result of Bedford-Klingenberg [BK] and Kru\v{z}ilin [K], there is a smooth 3-ball
$\mathcal B$ which is the disjoint union of holomorphic discs $\{D_\alpha \}$, such that  
$\partial {\mathcal B} = \tilde{S}$. By Chirka's theorem [C] we know that discs 
$D_\alpha$  are $C^\infty$-smooth to the boundary $\partial D_\alpha$  near the totally real part of
$\tilde{S}$  (i. e. outside of finitely many complex points of $\tilde{S}$) and, moreover, the boundary 
$\partial D_\alpha$  of each disc $D_\alpha$  is $C^\infty$-smooth at this part of $\tilde{S}$ and 
transversal there to the distribution $\{T_p^{\mathbb C}(\partial G)\}$ of complex tangencies to 
$\partial G$. Consider a disc $D_{\alpha_0}$  from the family $\{D_\alpha\}$ such that its boundary 
$\partial D_{\alpha_0}$ "touches" the curve $\Gamma '\,$ "for the first time" and let $p$ be a point of 
$\,\partial D_{\alpha_0} \cap \Gamma '$. More precisely, let $D_{\alpha_0} \subset {\mathcal B}\,$ be a
holomorphic disc with the property that $\partial D_{\alpha_0} \cap \Gamma' \neq \emptyset$ and such that 
for some connected component of the set $\,{\mathcal B} \backslash D_{\alpha_0}\,$, each holomorphic disc
$D_\alpha$, which is contained in this component, satisfies $\partial D_\alpha \cap \Gamma ' =\emptyset$. 
Now we can see that, on the one hand, since the curves $\partial D_{\alpha_0}$ and $\Gamma '$ are tangent to each
other at the point $p$, and since the curve $\Gamma '$ was chosen to be complex tangential, the curve 
$\partial D_{\alpha_0}$  is complex tangential at $p$. On the other hand, since the point $p$ is contained 
in the totally real part of $\tilde{S}$, the curve $\partial D_{\alpha_0}$  has to be transversal to 
$T_p^{\mathbb C}(\partial G)$. This gives the desired contradiction and completes the proof of the first part 
of the theorem.\\

\vskip0.15cm

\noindent
{\bf Remark}. In the special case when the boundary of the domain $G$ is diffeomorphic to a 3-dimensional sphere, the fact that the
closed complex tangential curve $\Gamma  \subset S\subset
\partial G$ mentioned above does not exist can also be deduced from the theorem
1 in [Be].
\vskip0.6cm

\noindent
{\large \bf {3. Proof of the second part of the theorem}}\\

\noindent
We prove the second part of the theorem in several steps. 
First, we construct a special arc $\,E \subset {\mathbb R}^2_{x,y}$ of positive 2-dimensional measure.
Then we define a function $G$ on $E$ such that $\,G \in C^{2-}(E)\,$ with the functions $G'_x (x,y) = y$ and $G'_y (x,y) = 0$ chosen to be the first 
derivatives of $G$ on $E$.
Next, following an idea of H. Whitney (see [Wh]), we construct a nonconstant function $H \in C^{2-}(E)$ with $H'_x (x, y) = 0$ and $H'_y (x,y) = 0$.
Using $G$ and $H$ we define a function $F \in C^{2-} (E)$ with $F'_x (x, y) = y$ and $F'_y (x, y) = 0$ which is zero at both endpoints of $E$.
Then, using $E$, we construct a Jordan curve $\tilde{E} \subset {\mathbb R}^2_{x, y}$ of positive 2-dimensional measure and, using $F$, we define a 
function $\tilde{F}$ on $\tilde{E}$ of class $C^{2-}(\tilde{E})$ with derivatives $\tilde{F}'_x (x,y) =y$ and
$\tilde{F}'_y (x,y) = 0$. Next, applying Whitney's extension theorem to the function $\tilde{F}$, we
contruct a 2-sphere $S^2 \subset {\mathbb R}^3_{x, y, z}$ of class $C^{2-}$ which contains  a Jordan curve of positive
2-dimensional measure such that at each point of this curve the tangent plane to $S^2$ coincides with the
corresponding  plane of the standard contact distribution in ${\mathbb
  R}^3_{x, y, z}$. Finally, using the Darboux
theorem, we embed this sphere into the boundary $bG$  of the given strictly pseudoconvex domain $G$.

\vskip0.5cm
\noindent
{\bf 1. Construction of the arc $\bf E.$\,\,}  First, we define for each $\alpha  \in (0,1)$ an auxiliary set
$$ {\mathbb E}^\alpha = \bigg( \bigg[ 0, \frac{1-\alpha }{2} \bigg] \cup \bigg[ \frac{1+\alpha }{2}, 1 \bigg] \bigg) \times 
\bigg( \bigg[ 0, \frac{1-\alpha }{2} \bigg] \cup \bigg[ \frac{1+\alpha }{2}, 1 \bigg] \bigg) \cup $$\
$$\bigg( \{0\} \times [ 0, 1 ]\bigg)  \cup \bigg( [ 0, 1 ] 
\times \bigg\{ \frac{1+\alpha }{2} \bigg\} \bigg) \cup 
\bigg( \{1\} \times [ 0, 1 ] \bigg).$$

\vskip0.5cm 

We denote $A = (0, 0), B = (1, 0), Q_0 = [0, \frac{1 - \alpha }{2}] \times [0, \frac{1 - \alpha }{2}], Q_1 = [0, \frac{1 - \alpha }{2}] \times
[\frac{1 +\alpha }{2}, 1], Q_2 = [\frac{1 + \alpha }{2}, 1] \times [\frac{1 + \alpha }{2}, 1]$, and $Q_3 = [\frac{1 + \alpha }{2}, 1] \times 
[ 0, \frac{1 - \alpha}{2}]$.
Further, we denote $A_0 = A = (0, 0), B_0 = (0, \frac{1 - \alpha }{2}), A_1 = (0, \frac{1 + \alpha }{2}), B_1 = (\frac{1 - \alpha}{2}, 
\frac{1 + \alpha}{2}), A_2 = (\frac{1 + \alpha}{2}, \frac{1 + \alpha}{2}), B_2 = (1, \frac{1 + \alpha}{2}), A_3 = (1, \frac{1 - \alpha}{2})$, and 
$B_3 = B = (1, 0)$ (see the set ${\mathbb E}^\alpha$ in Fig. 1).
\begin{figure}[ht]
\begin{center}\input{fig1.tex}\end{center}
\caption{The set $\mathbb{E}^\alpha$.}
\end{figure}
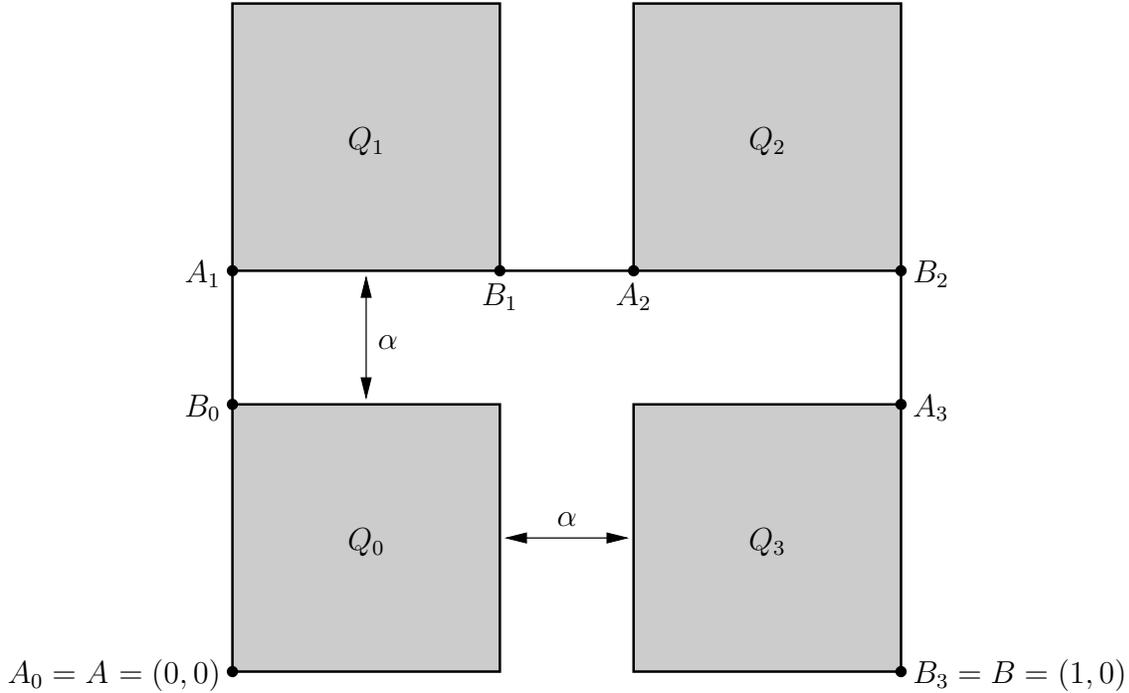

\vskip0.2cm 

To define the set $E$ we consider the sequence $\alpha_n = \frac{1}{(n+1)^2}, n = 1, 2, \ldots $.
We construct the set $E$ as the intersection of a decreasing sequence of compact sets $E_n$ which will be defined inductively. 
We set $E_1 = {\mathbb E}^{\alpha_1}$.
To define the set $E_2$ we consider the image $\tilde{{\mathbb E}}^{\alpha_2}$ of the set  ${\mathbb E}^{\alpha_2}$ under the homothety with 
coefficient $\frac{1 - \alpha_1}{2}$.
Then for each $i = 0, 1, 2, 3 \,$ we consider the set ${\mathbb E}_i$ obtained from the set  $\tilde{{\mathbb E}}^{\alpha_2}$ by an orthogonal 
transformation (if necessary) and translation in such a way that the image of the points $A$ and $B$ will coincide with the points $A_i$ and $B_i$, 
respectively.
It is easy to see that we need an orthogonal transformation  only for $ i = 0, 3$. The set $E_2$ is obtained from the set $E_1$ by replacing
for each $i = 0, 1, 2,
3 \,$ the square $Q_i$ by the corresponding set ${\mathbb E}_i$. For each $j = 0, 1, 2, 3 \,$ we denote by $Q_{ij},
A_{ij}$ and $B_{ij}$ the images of the square $Q_j$ and the points $A_j, B_j$ in the corresponding  set
${\mathbb E}_i$, respectively (the set $E_2$ is shown in Fig. 2).
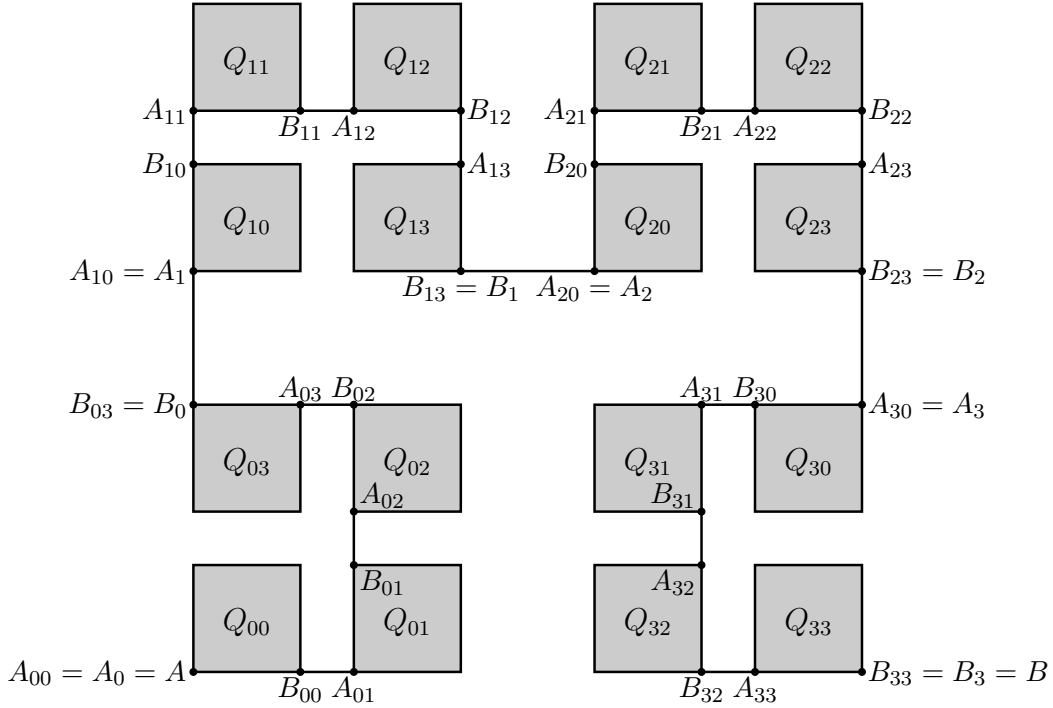
\begin{figure}[ht]
\begin{center}\input{fig2.tex}\end{center}
\caption{The set $E_2$.}
\end{figure}

\vskip0.2cm

To describe the inductive step of our construction we assume that the set $E_n$ is already constructed and define the set $E_{n+1}$.
Consider the image $\tilde{{\mathbb E}}^{\alpha_{n + 1}}$ of the set ${\mathbb E}^{\alpha_{n + 1}}$ under the
homothety with coefficient 
$ \prod\limits^n_{i = 1} \big( \frac{1 - \alpha_i }{2} \big)$.
Then for each multiindex $(i_1, i_2, \ldots , i_n), i_j = 0, 1, 2, 3, \, j =
1, 2, \ldots , n $, consider the set ${\mathbb E}_{i_1 \ldots i_n}$ 
obtained from the set $\tilde{\mathbb E}^{\alpha_{n+1}}$ by an orthogonal transformation (if necessary) and translation in such a way that the image 
of the points $A$ and $B$ will coincide with the points $A_{i_1 \ldots i_n}$ and $B_{i_1 \ldots i_n}$, respectively.
The set $E_{n+1}$ is obtained from the set $E_n$ by replacing each square $Q_{i_1 \ldots i_n}$ by the corresponding set ${\mathbb E}_{i_1 \ldots i_n}$.
For each $i_j = 0, 1, 2, 3, \, j = 1, 2, \ldots , n+1$, we denote by  $Q_{i_1 \ldots i_{n+1}}, A_{i_1 \ldots i_{n+1}}$ and $ B_{i_1 \ldots i_{n+1}}$ the 
images of the square $Q_{i_{n+1}}$ and the points $A_{i_{n+1}}, B_{i_{n+1}}$ in the corresponding set ${\mathbb E}_{i_1 \ldots i_n}$, respectively.
Note, that for each multiindex $(i_1, \ldots , i_n)$ one has $A_{i_1 \ldots i_n 0} = A_{i_1 \ldots i_n}$ and $B_{i_1 \ldots i_n
3} = B_{i_1 \ldots i_n}$.

\vskip0.2cm

Since $\{ E_n \}$ is a decreasing sequence of compact sets, $E = \bigcap\limits^\infty_{n = 1} E_n$ is a nonempty compact subset of ${\mathbb R}^2_{x, y}$.
It is easy to see that it is an arc (see the set $E$ in Fig. 3).
\begin{figure}[ht]
\begin{center}\input{fig3.tex}\end{center}
\caption{The set $E$.}
\end{figure}
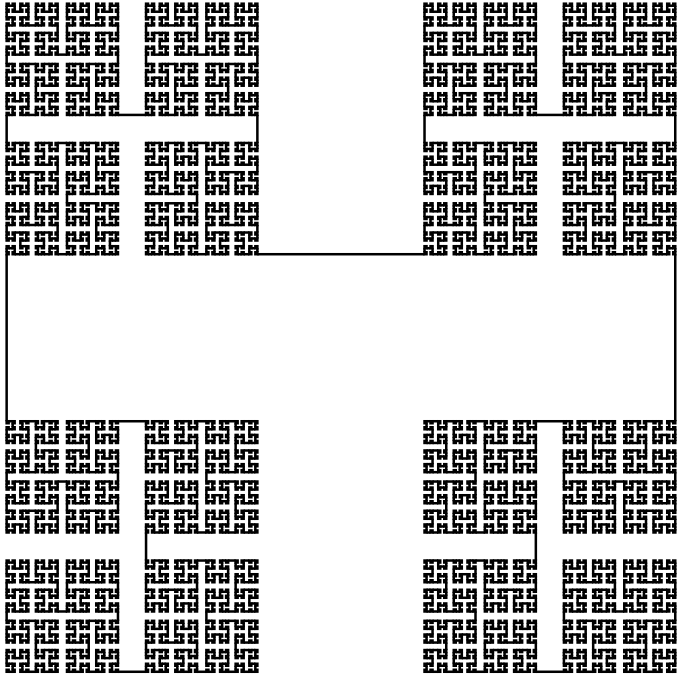

\vskip0.2cm

To estimate the area of the set $E$ we observe that 
\begin{eqnarray} 
{\rm Area} \, \, (E_n) = (1 - \alpha_n)^2 \,\, {\rm Area} \, \, (E_{n-1}) & = & \prod\limits^n_{k = 1} (1 - \alpha_k)^2  = \prod\limits^n_{k = 1} 
\bigg( 1 - \frac{1}{(k + 1)^2} \bigg)^2 \nonumber\\ & = & \bigg( \frac{1}{2} \bigg( 1 + \frac{1}{n+1} \bigg) \bigg)^2 > \frac{1}{4} \nonumber  
\end{eqnarray} 
for every $n = 1, 2,  \ldots$.
Hence, the set $E$ has a positive 2-dimensional measure (and, moreover, ${\rm Area} \, \, (E) = \frac{1}{4}$).

\vskip0.5cm
\noindent
{\bf 2. Definition and properties of the function $\bf G.$\,\,} For each $n = 1, 2, \ldots $ let $\Omega_n$ be the connected component of the set 
$(0, 1) \times (-1, 1) \setminus E_n$ containing the square 
$(0, 1) \times (-1, 0)$ and let $J_n = \partial \Omega_n \cap ([0,1] \times [0,1])$.
On each curve $J_n$ we define a function $G_n$ in the following way.
For a point $p \in J_n$ we denote by $J^p_n$ a part of $J_n$ with initial point $A$ and endpoint $p$ and then we set 
$G_n (p) = \int \limits_{J^p_n}  y dx$.

\vskip0.2cm
We will need the following estimate on the function $G_n$.

\vskip0.2cm
\noindent
{\bf Lemma}. Let $Q_{i_1 i_2 \ldots i_m}$ be a square obtained after $m \geq 6$ steps of our construction above. Then for each 
natural number $n > m$ and all points $p, q \in J_n \cap  Q_{i_1 i_2 \ldots i_m}$ one has 
\begin{equation} \mid G_n (q) - G_n (p) - \int_{[p, q]} ydx \mid \,\, < Area (Q_{i_1 i_2 \ldots i_m}) = \frac{1}{2^{2m+2}} 
\bigg(\frac{m+2}{m+1}\bigg)^2 \,\, .
\end{equation}

\vskip0.2cm
\noindent
{\bf Proof}. We proceed by induction on $m$ having the number $n$ fixed and decreasing the number $m$. For each pair of points
$r, s \in J_n$ we denote by $J_n^{r,s}$ the part of the curve $J_n$ with
initial point $r$ and endpoint $s$. For $m = n-1$ we observe from the 
definition of the function $G_n$ and Green's theorem that $G_n (q) - G_n (p) - \int_{[p,q]} ydx$ represents the sum (with signes) of the areas 
of domains bounded by the arc $J_n^{p,q}$ and the segment $[p, q]$ (see Fig. 4).

\vskip0.5cm

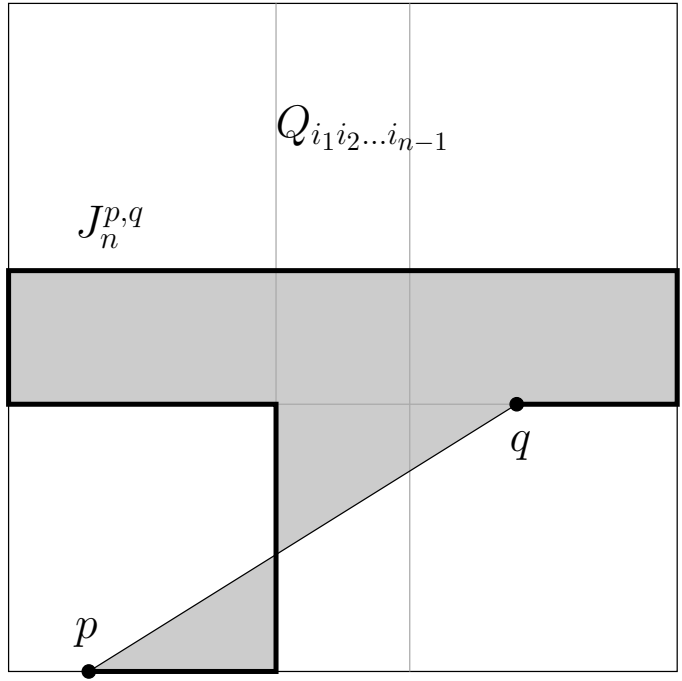
\begin{figure}[ht]
\begin{center}\input{fig5.tex}\end{center}
\caption{The curve $J_n^{p,q}$ in the square $Q_{i_1 i_2 \ldots i_{n-1}}$.}
\end{figure}

\vskip0.5cm

Since these domains are contained in the square $Q_{i_1 i_2 \ldots i_{n-1}}$, we conclude that the inequality (1) holds true in this case. 
The same argument can not be applied in the case of general $m$ due to the fact that the corresponding domains will not be disjoint 
anymore and, therefore, their areas will be counted with multiplicities. That is why in the case of an arbitrary $m$ we need to improve 
the argument above and to use the inductive procedure.

\vskip0.2cm

Assume that the inequality (1) holds true for all squares $Q_{i_1 i_2 \ldots i_m}$ and all points 
$p,q \in J_n \cap Q_{i_1 i_2 \ldots i_m}$ with $m \geq m_0 +1$ ($n$ being fixed). Consider a square $Q_{i_1 i_2 \ldots i_{m_0}}$ and points 
$p,q \in J_n \cap Q_{i_1 i_2 \ldots i_{m_0}}$. Observe first that the segment $[p, q]$ can intersect at most three of the squares 
$Q_{i_1i_2 \ldots i_{m_0} i_{m_0 +1}}, \,\, i_{m_0 +1} = 0, 1, 2, 3$. We assume in what follows that it intersects exactly three 
of them (in the other cases the argument is easier, since some terms in our estimates will disappear), namely the squares 
$Q_{i_1 i_2 \ldots i_{m_0}0}, Q_{i_1 i_2 \ldots i_{m_0}1}$ and $Q_{i_1 i_2 \ldots i_{m_0}2}$ (the cases of another three squares 
can be treated similarly). Consider an orientation of the segment $[p, q]$
from the point $p$ to the point $q$ and for each $i = 0, 1, 2$ 
denote by $p_i$ and $p'_i$ the first and the last, respectively, (according to the orientation of $[p, q]$) points of the set 
$[p, q] \cap  J^{p,q}_n \cap Q_{i_1 i_2 \ldots i_{m_0}i}$. Note, that $p_0 =
p$ and $p'_2 =q$. For each $i = 0, 1, 2$ we have by 
our inductive hypotesis that
\begin{equation} \mid G_n (p'_i) - G_n(p_i)- \int_{[p_i,\,  p'_i]} ydx \mid
  \,\, < Area (Q_{i_1 i_2 \ldots i_{m_0}i}) < \frac{1}{4} \,  
Area (Q_{i_1 i_2 \ldots i_{m_0}}) \,\, .
\end{equation} 

Since the segment $[p'_0, p_1]$ intersects the curve $J_n$ only at its endpoints, we conclude from the Green's theorem that 
$G_n (p_1)- G_n (p'_0) - \int_{[p'_0, \,\, p_1]} ydx$
represents the area of the domain $\Sigma _1$ bounded by the curve
$J_n^{p'_0,\, p_1}$ and the segment $[p'_0, p_1]$ (see Fig. 5).

\vskip0.5cm

\begin{figure}[ht]
\begin{center}\input{fig6b.tex}\end{center}
\caption{The sets $\Sigma _1$ and $\Omega _n^{i_1 i_2 \ldots i_{m_0}}$ in the spuare $Q_{i_1 i_2 \ldots i_{m_0} }$.}
\end{figure}
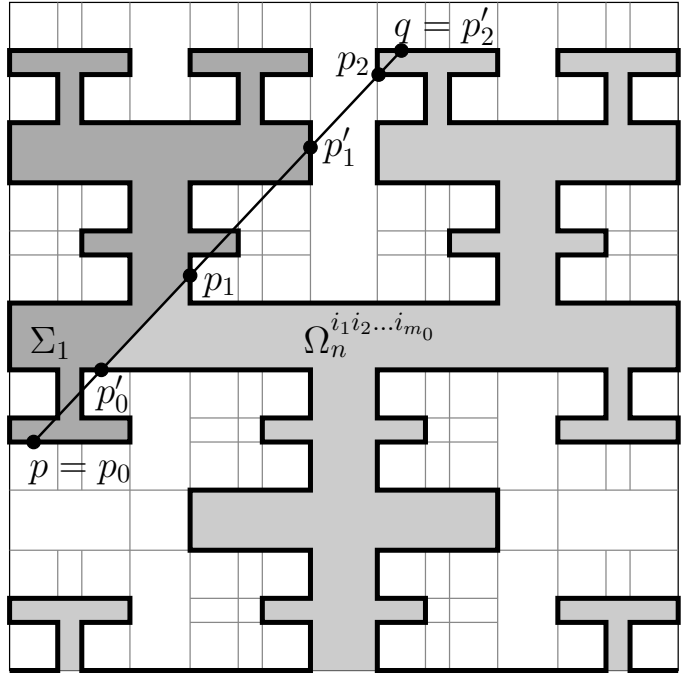

\vskip0.5cm

Since the domain $\Sigma _1$ is contained in the set $\Omega _n^{i_1 i_2 \ldots i_{m_0}} = \Omega _n \cap  Q_{i_1 i_2 \ldots i_{m_0}}$, 
and since $\Omega _n \cap  Q_{i_1 i_2 \ldots i_{m_0}} \subset  Q_{i_1 i_2 \ldots i_{m_0}} \setminus\ E_n$, it follows that
\begin{equation}  0 < G_n (p_1) - G_n (p'_0) - \int_{[p'_0, \,\, p_1]} ydx = Area (\Sigma _1) < Area (\Omega _n \cap Q_{i_1 i_2 \ldots i_{m_0}}) $$
$$ < Area (Q_{i_1 i_2 \ldots i_{m_0}}) - Area (E_n\cap  Q_{i_1, i_2 \ldots i_{m_0}}) = Area (Q_{i_1 i_2 \ldots i_{m_0}}) (1 - (1 - \alpha _{m_0 +1})^2 \ldots (1 - \alpha_n)^2) $$
$$ =Area (Q_{i_1 i_2 \ldots i_{m_0}}) \bigg(1 - \bigg(1 -
\frac{1}{(m_0  +2)^2}\bigg)^2 \ldots \bigg(1 - \frac {1}{(n+1)^2}\bigg)^2
\bigg) $$
$$ = Area (Q_{i_1 i_2 \ldots i_{m_0}})\bigg(1 - \bigg(\frac{(m_0 + 1) (n+2)}{(m_0  +2)(n+1)}\bigg)^2\bigg) < Area (Q_{i_1 i_2 \ldots i_{m_0}}) 
\bigg(1 - \bigg(\frac{m_0 +1}{m_0 +2}\bigg)^2\bigg) $$
$$<\frac {2}{m_0 + 2} \,\, Area (Q_{i_1 i_2 \ldots i_{m_0}})\,\, . 
\end{equation} 
Similarly, $G_n (p_2) - G_n (p'_1) - \int_{[p'_1, \,\, p_2]} ydx$ represents the area with the negative sign (due to the orientation of its boundary) of 
the domain $\Sigma _2$ bounded by the curve $J_n^{p'_1, \, p_2}$ and the segment $[p'_1, p_2]$. Let $\Pi$ be the smallest rectangle in 
$Q_{i_1 i_2 \ldots i_{m_0}}$ that coutains the squares $Q_{i_1 i_2 \ldots i_{m_0}12}, Q_{i_1 i_2 \ldots i_{m_0}13}, Q_{i_1 i_2 \ldots i_{m_0}20}$ 
and $Q_{i_1 i_2 \ldots i_{m_0}21}$ \,\,(see Fig. 6).

\vskip0.5cm

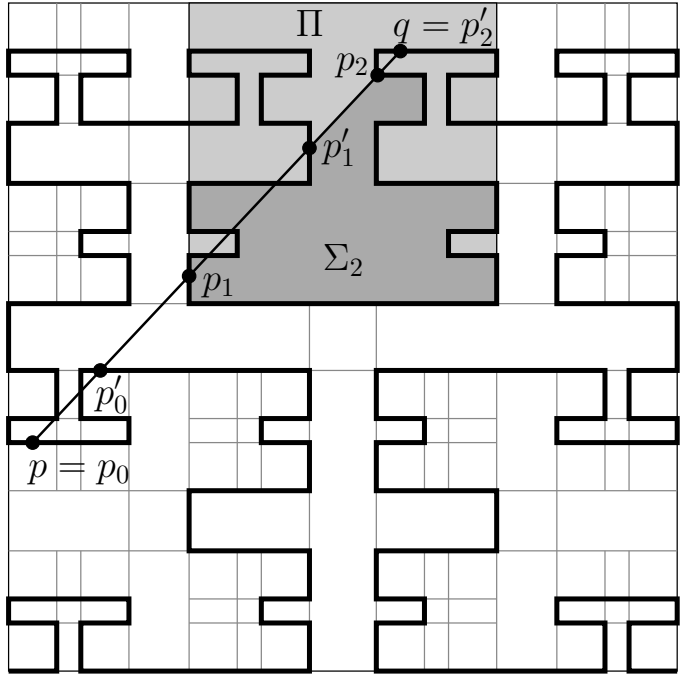
\begin{figure}[ht]
\begin{center}\input{fig7b.tex}\end{center}
\caption{The sets $\Sigma _2$ and $\Pi $ in the square $Q_{i_1 i_2 \ldots i_{m_0}}$.}
\end{figure}

\vskip0.5cm
Since the domain $\Sigma _2$ is a subdomain of $\Pi $, and since we have an obvious estimate 
$Area (\Pi) < \frac{1}{4} \, Area (Q_{i_1 i_2 \ldots i_{m_0}})$, we conclude that
%
%
\begin{equation} 0> G_n(p_2) - G_n (p'_1) - \int_{[p'_1, \,\, p_2]} ydx = - Area (\Sigma _2) > - Area (\Pi )  $$
$$ > - \frac{1}{4} Area (Q_{i_1 i_2 \ldots i_{m_0}}) \,\, .
\end{equation}  

\vskip0.2cm

\noindent
Taking into account different signs of the terms in the inequalities (3) and (4) we have for $m_0  \geq  6$ the following estimate
\begin{equation}  \mid (G_n (p_1) - G_n (p'_0) - \int _{[p'_0, \,\, p_1]}
  ydx)+ (G_n (p_2) - G_n (p'_1) - \int_{[p'_1, \,\, p_2]} ydx ) \mid  $$ 
$$< \,\, \frac{1}{4} \, Area (Q_{i_1 i_2 \ldots i_{m_0}}) \,\, .
\end{equation} 

\vskip0.2cm
\noindent
Finally, the inequalites (2) and (5) imply that 
$$\mid G_n (q) - G_n (p) - \int _{[p,q]} ydx \mid \, < \, \sum^2_{i=0} \mid (G_n
(p'_i) - G_n (p_i) - \int_{[p_i, \,\, p'_i]} ydx \mid  +$$
$$+ \mid (G_n (p_1) - G_n (p'_0) - \int _{[p'_0, \,\, p_1]} ydx)+ (G_n (p_2) -
G_n (p'_1) - \int_{[p'_1, \,\, p_2]} ydx \mid  $$ 
$$<\,\, \frac{3}{4} \, Area (Q_{i_1 i_2 \ldots i_{m_0}}) + \frac{1}{4} \,
Area (Q_{i_1 i_2 \ldots i_{m_0}}) = Area (Q_{i_1 i_2 \ldots i_{m_0
  }}) \,\, .$$
This completes the proof of the lemma.\bx\\


Now we define the function $G$ on the arc $E$.
Let $p$ be a point of $E$.
Consider a sequence of points $p_n \in J_n,\, n = 1,2, \ldots $, such that $p_n \rightarrow p$, as $n \rightarrow  \infty $ and set 
$G(p) =  \lim\limits_{n \rightarrow  \infty } G_n (p_n)$.
\vskip0.2cm

To prove that the function $G$ is well defined we consider the Cantor set ${\mathbb Q} \overset{def}{=}  \bigcap\limits_{n=1}^\infty  
\bigcup\limits_{(i_1, i_2, \ldots, i_n)}  Q_{i_1 i_2 \ldots i_n}$ and observe that since the set $E \setminus\ {\mathbb Q}$ is 
constituted by horizontal and vertical segments, and in view of the 
definition of the functions $G_n$, it is enough to prove that the sequence $G_n (p_n)$ converges for $p_n \rightarrow p$ with $p, p_n \in {\mathbb Q}, 
\,  n = 1, 2, \ldots$. Let $p_{n_1}$ and $p_{n_2}$, $ n_1 > n_2$, be two points of
our sequence $p_n \in J_n \bigcap {\mathbb Q}, \, n = 1, 2, \ldots$, and let 
$Q_{i_1 i_2 \ldots i_m}$ be the smallest of the described above squares that contains both of these points. Observe that the vertex $A_{i_1 i_2 \ldots i_m}$ 
of the square $Q_{i_1 i_2 \ldots i_m}$ is contained in both curves $J_{n_1}$ and $J_{n_2}$. It follows then from Green's theorem that $G_{n_1} 
(A_{i_1 i_2 \ldots i_m}) - G_{n_2} (A_{i_1 i_2 \ldots i_m})$ represents the area of the set bounded by the curves $J_{n_1}^{A, A_{i_1 i_2 \ldots i_m}}$ and 
$J_{n_2}^{A, A_{i_1 i_2 \ldots i_m}}$. Since this set is a subset of the set $ \Omega _{n_1} \setminus\ \bar \Omega _{n_2}$, we conclude that
\begin{equation} \mid G_{n_1} (A_{i_1 i_2 \ldots i_m}) - G_{n_2} (A_{i_1 i_2
    \ldots i_m})\mid \,\, < \, Area ( \Omega _{n_1} \setminus\ \bar \Omega
    _{n_2}) \,\, .
\end{equation}  
\noindent
For each $ j= 1, 2$ we have (by the lemma):
$$ \mid G_{n_j} (p_{n_j}) - G_{n_j} (A_{i_1 i_2 \ldots i_m}) - \int_{[A_{i_1
    i_2 \ldots i_m}, \,\, p_{n_j}]} ydx \mid \,\, < \, Area ( Q_{i_1 i_2
    \ldots i_m}) \,\, .$$
\noindent
If we denote by $l_m$ the length of a side of $Q_{i_1 i_2 \ldots i_m}$ (it depends only on $m$), then the last inequality gives us 
\begin{equation}  \sum\limits_{j = 1}^{2} \mid G_{n_j} (p_{n_j}) - G_{n_j}
  (A_{i_1 i_2 \ldots i_m}) \mid \, < \, \sum\limits_{j = 1}^{2} \mid \int_{[A_{i_1 i_2 \ldots i_m}, 
\,\, p_{n_j}]}ydx \mid + 2 Area (Q_{i_1 i_2 \ldots i_m}) $$
$$< 2 \, l_m + 2 Area (Q_{i_1 i_2 \ldots i_m})\,  < \, 3 \, l_m = 3 \prod\limits^m_{i = 1}
\bigg( \frac{1-\alpha_i}{2} \bigg) = \frac{3}{2^{m+1}} \cdot \frac{m+2}{m+1}
\,\, .
\end{equation} 
\noindent
It follows now from (6) and (7) that

$$\mid G_{n_1} (p_{n_1}) - G_{n_2} (p_{n_2})  \mid \,\, \leq \,\, \sum\limits_{j = 1}^{2} \mid G_{n_j} (p_{n_j}) - G_{n_j} (A_{i_1 i_2 \ldots i_m}) \mid + $$
$$ \mid  G_{n_1} (A_{i_1 i_2 \ldots i_m}) -  G_{n_2} (A_{i_1 i_2 \ldots i_m})
\mid \,\, < \,\, \frac{3}{2^{m+1}} \cdot \frac{m+2}{m+1} 
+ Area (\Omega _{n_1} \setminus\ \bar\Omega _{n_2})\,\, .$$

\noindent
Since for points $p, p_{n_1}, p_{n_2} \in {\mathbb Q}$ the number $m$ will tend to infinity as $p_{n_1}, p_{n_2} \rightarrow  p$, and since 
$Area (\Omega _{n_1} \setminus\ \bar\Omega _{n_2})\rightarrow 0$ as $n_1, n_2 \rightarrow \infty $, we conclude that the limit of the sequence $G_n (p_n)$ 
as $n \rightarrow  \infty $ exists and, therefore, the function $G$ is well defined.

\vskip0.2cm

Next, we prove that $G \in C^{2-} (E)$ with $G'_x (x,y) = y$ and $G'_y (x,y) = 0$ chosen to be the first derivatives of $G$ at
each point $(x, y) \in E$. To do this rigorously we recall the definition of a function belonging to the class $C^{2-} (E)$
(this definition is due to H. Whitney.  Further details can be found, for example, in [S]).

\vskip0.5cm
\noindent
{\bf Definition.}
{\it Let $E$ be a compact subset of $\,\,{\mathbb R}^2_{x,y}$ and let $f$ be a function defined on $E$.
We say that $f$ belongs to the class $C^{2-} (E)$ if there exist bounded
functions $f'_x$ and $f'_y$ defined on $E$  with the property that for each $  \varepsilon > 0$ 
there is a constant $M$ such that 
\begin{equation}
|f (x + \Delta x, y + \Delta y) - f(x,y) -f'_x (x, y) \Delta  x - f'_y (x, y) \Delta  y | 
\le M (| \Delta  x| + | \Delta  y|)^{2-\varepsilon} 
\label{lg8} \end{equation}
\vskip0.2cm
\noindent
for all $\,\,(x, y), (x + \Delta x, y + \Delta  y) \in E$.}

\vskip0.5cm

To prove that $G \in C^{2-} (E)$ we consider two points $\,\,p, p + \Delta p \in E$. 
Since the function $G$ obviously is smooth on each of the segments $[B_{i_1 \ldots i_n}, A_{i_1 \ldots (i_n + 1)}], i_n = 0, 1, 2,$  
and satisfies condition (8) with $G'_x (x, y) = y$ and
$G'_y (x, y) = 0$ there, the general case, when $p, p + \triangle p \in E$, can be easily reduced to the case when $p$ and $p + \triangle p$ 
are contained in the Cantor set ${\mathbb Q}$. That is why we assume in what follows that $p, p + \triangle p \in {\mathbb Q}$. 
Consider as above a number
$m$ such that $p, p + \Delta p$ belong to a square $Q_{i_1 \ldots i_m}$ for some
indices $(i_1, \ldots , i_m)$, but not to a smaller square $Q_{i_1, \ldots i_m i_{m+1}}, i_{m+1} = 0, 1, 2, 3$.
Since $p$ and $p + \Delta p$ belong to different squares $Q_{i_1 \ldots i_m i_{m+1}}$ and $Q_{i_1 \ldots i_m i'_{m+1}}$,
it follows that the distance  between these points is not less than the minimal distance between $Q_{i_1 \ldots i_m i_{m+1}}$
and $Q_{i_1 \ldots i_m i'_{m+1}}$, that is,
\begin{eqnarray}
|\Delta  p| \ge \alpha_{m+1} \left(\frac{1- \alpha_1}{2}\right) \ldots \left(\frac{1- \alpha_m}{2}\right) &=& \frac{1}{(m+2)^2} \cdot \frac{1}{2^m} 
\prod\limits^m_{k = 1} \bigg( 1 - \frac{1}{(k + 1)^2} \bigg)  \nonumber\\ 
&=& \frac{1}{2^{m + 1}} \cdot \frac{1}{(m + 1)(m+2)} \,\, . \label{lg9} 
\end{eqnarray}

\vskip0.5cm

Now we estimate the left hand side of the condition (8) for our function $G$

$$ {\mathcal L}_G (p, p + \Delta  p) \overset{def}{=} G (p + \Delta  p) - G (p) - G'_x (p) \Delta x - G'_y (p) \Delta y = 
G (p + \Delta  p) - G (p) - y \Delta  x, $$

\vskip0.2cm
\noindent
where $p = (x, y)$ and $\,\Delta  p = (\Delta x, \Delta y)$. 
It is easy to see that 
$$\int_{[p, \,\, p + \Delta p]}  y dx = y \Delta  x + \frac{1}{2} \Delta  x \Delta y, $$ 
hence
\begin{eqnarray}
 | {\mathcal L}_G (p, p + \Delta  p)| \, \le \, | G (p + \Delta p) - G (p) -
 \int_{[p, \,\, p + \Delta p]} y dx | + \frac{1}{2}\, |\Delta x||\Delta
y| \,\, . \label{lg10} 
\end{eqnarray}
\vskip0.5cm

It follows from the lemma above that 
\begin{eqnarray} 
{\Big|G (p + \Delta  p) - G (p) - \int_{[p, \,\, p + \Delta p]} y dx \Big| \le  } 
 \,\, \frac{1}{2^{2m+2}} \left( \frac{m + 2}{m + 1} \right)^2 \,\, .\label{lg11} 
\end{eqnarray}
\vskip0.5cm
Since $\,\,p, p + \Delta p \in Q_{i_1 \ldots i_m}$, one can estimate $|\Delta  x|$ and $|\Delta  y|$ from above by the length 
of the side of $Q_{i_1 \ldots i_m}$, that is, by $\frac{1}{2^{m + 1}} \left( \frac{m + 2}{m + 1}\right)$. 
Therefore, we have by (10) and (11) that 
\begin{equation}| {\mathcal L}_G (p, p + \Delta  p) | \, \le \,\, \frac{3}{2}  \cdot \frac{1}{2^{2m + 2}}  \left(\frac{m + 2}{m + 1}
\right)^2 \,\, .\label{lg12}
\end{equation}

\vskip0.5cm

Finally, we conclude from the estimates (9) and (12) that in order to prove, that $G$ satisfies condition (8)
we only need to verify that for each 
$\varepsilon > 0 $ there is a constant $M$ such that
$$ \frac{3}{2}  \cdot \frac{1}{2^{2m + 2}} \left( \frac{m+2}{m+1} \right)^2 \le  M \left( \frac{1}{2^{m + 1}}  \cdot \frac{1}{(m + 1)(m+2)} 
\right)^{2-\varepsilon}  \,\,{\rm as} \,\, m \rightarrow \infty $$
which is equivalent to the inequality
$$\frac{1}{(2^\varepsilon)^{m + 1}} \le M \, \frac{2}{3} \left( \frac{1}{(m+1)(m+2)} \right)^{2-\varepsilon} \left( \frac{m + 1}{m + 2} \right)^2  \,\,
{\rm as} \,\, m \rightarrow  \infty.$$
The last inequality is obviously satisfied, since the left hand side tends to zero much faster than the right hand side, as $m \rightarrow  \infty $. 
This proves that the function $G$ belongs to the class $C^{2-} (E)$.

\vskip0.5cm
\noindent
{\bf 3. Definition and properties of the functions $\bf H$ and $\bf F$.\,\,} First, we define the function $H$ on the Cantor set ${\mathbb Q}$. 
Each point $p$ in this set is uniquely determined as the intersection of the decreasing sequence $Q_{i_1} \supset Q_{i_1 i_2} \supset  
Q_{i_1 i_2 i_3} \supset \ldots$ of the squares $Q_{i_1 \ldots i_n}$.
Then, we define the value of $H$ at the point $p$ as $H(p) =  \sum\limits_{n = 1}^{\infty } \frac{i_n}{4^n}$.
It is easy to see that for each $i_n = 0, 1, 2$ one has $H ( A_{i_1 \ldots i_{n-1} (i_n +1)}) =  \sum\limits_{k = 1}^{n} \frac{i_k}{4^k} + \frac{1}{4^n}$ 
and $H (B_{i_1 \ldots i_{n-1} i_n}) =  \sum\limits_{k = 1}^{n} \frac{i_k}{4^k} +  \sum\limits_{k = n + 1}^{\infty }
\frac{3}{4^k}  = \sum\limits_{k = 1}^{n} \frac{i_k}{4^k} + \frac{1}{4^n} $, therefore, we can extend the function $H$ as a
constant to each segment 
$[B_{i_1 \ldots i_{n-1} i_n}, A_{i_1 \ldots i_{n - 1} (i_n + 1)} ], i_n = 0, 1, 2 $, with the value  
$\sum\limits_{k = 1}^{n} \frac{i_k}{4^k} + \frac{1}{4^n}$ there.
This defines the function $H$ on the whole set $E$.

\vskip0.2cm

Now we show that $H \in C^{2-} (E)$ with the functions $H'_x (x, y) = 0$ and $H'_y (x, y) = 0$ chosen to be the first derivatives of $H$ on $E$. 
We proceed in the same way as in the case of the function $G$, namely, we consider two points $p, p + \Delta p \in E$.
Since, by definition, $H$ is a constant on each of the intervals constituting the set $E \setminus {\mathbb Q}$, we only need to verify that the 
function $H$ satisfies condition (8) for points $\,p, p + \Delta p \in {\mathbb Q}$.
Let, as above, $m$ be a number such that $\,p, p + \Delta p \in Q_{i_1 \ldots i_m}$, but $\, p, p + \Delta p \not\in Q_{i_1 \ldots
i_m i_{m + 1}}$  for any $i_{m + 1} = 0, 1, 2, 3$.
Then  the definition of $H$ gives us that $|H (p + \Delta p) - H (p) | \le \frac{1}{4^m}$.
Hence, by estimate (9) , it is enough to show that for each $\varepsilon > 0$ there is $M$ such that
$$\frac{1}{4^m} \le M \left( \frac{1}{2^{m + 1}} \cdot  \frac{1}{(m + 1)(m+2)} \right)^{2 - \varepsilon} \,\,{\rm as} \,\, m
\rightarrow  \infty , $$ which is obviously true with the same argument as above for the function $G$.

\vskip0.2cm

To define the function $F$ on the set $E$ we first note that by definition of $H$ one has $H (A) = 0$ and $H (B) = 1$.
Then, since by definition of $G$ we have $G (A) = 0$, there is a constant $C$ such that for the function $F = G + CH$ one has $F(A) = 0$ and $F (B) = 0$.
Finally, we observe that since $G \in C^{2-}(E)$ with $G'_x (x, y) = y$ and $G'_y (x, y) = 0$, and since $ H \in C^{2-}(E)$ with $H'_x (x, y) = 0$ and 
$H'_y (x, y) =0$, it follows that $F \in C^{2-} (E)$ with $F'_x (x, y) = y$ and $F'_y (x, y) = 0$ at each point $(x, y) \in E$.

\begin{figure}[ht]
\begin{center}\input{fig4a.tex}\end{center}
\caption{The set $\tilde E$.}
\label{fig4}
\end{figure}
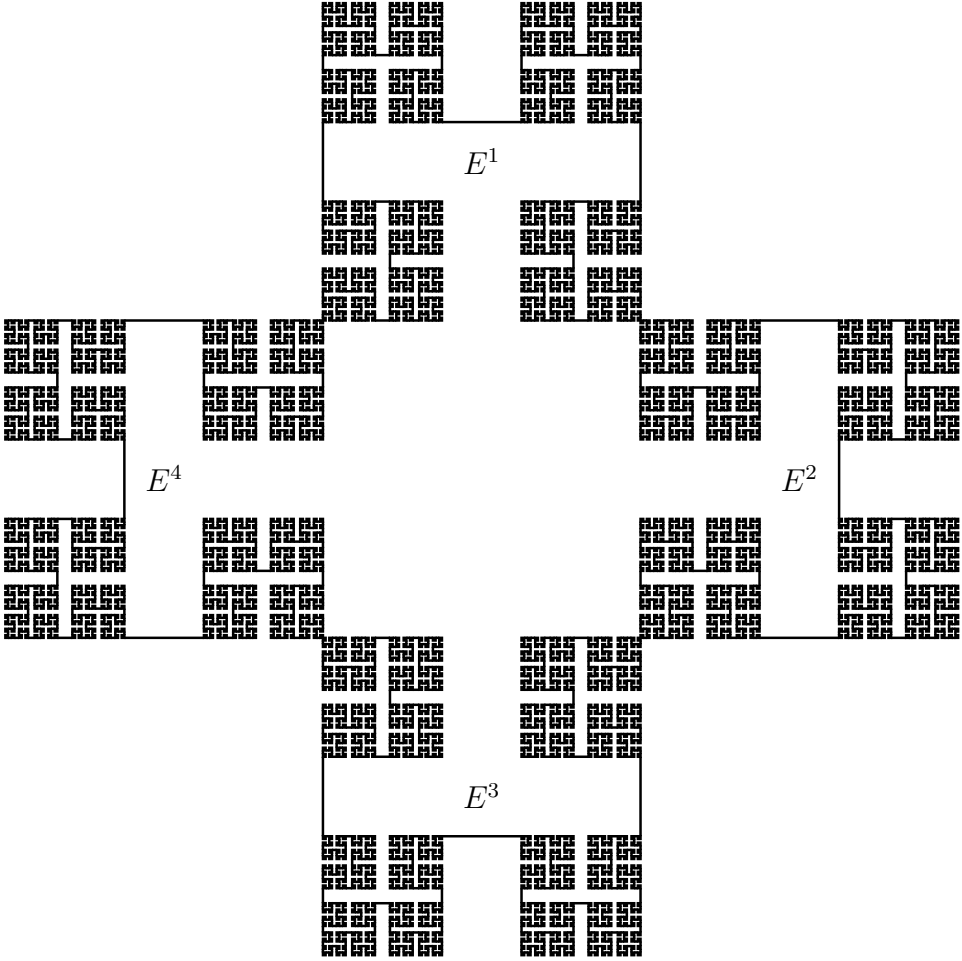

\vskip0.5cm
\noindent
{\bf 4. Construction of the sphere  $\bf {S \subset \partial G}$.\,\,} Let ${\mathbb A}$ be the linear transformation of $\,{\mathbb R}^2_{x, y}$ 
represented by the matrix $\displaystyle {\,\,\,\,0 \,\,\,\, 1\, \choose -1 \,\,\,\,
0\,} $.  Consider the sets $E^1 = E + \vec{e}_y, E^2 = -{\mathbb A} E   + \vec{e}_x + \vec{e}_y, E^3 = -E + \vec{e}_x$ and $E^4
= {\mathbb A} E$, where 
$\vec{e}_x$ and $\vec{e}_y$ are the unit vectors in the coordinate directions $x$ and $y$, respectively, and then define the
set $\tilde{E}$ as 
$\tilde{E} = \bigcup\limits_{i=1}^4 E^i$ (the set $\tilde{E}$ is shown
in Fig. \ref{fig4}). 

It is easy to see that $\tilde{E}$ is a Jordan curve of positive 2-dimensional measure in ${\mathbb R}^2_{x, y}$.
Applying to each of the sets $E^i, i = 1, 2, 3, 4$, a construction similar to the one that we had above for the function $F$ on the set $E$, 
we will get functions $F^i$ defined on the corresponding sets $E^i$ with the properties:

\vskip0.5cm

1) $F^i \in C^{2-} (E^i)$, 
\vskip0.5cm
2) $ \frac{\partial  F^i}{\partial  x} (x, y) = y\,\,$ and $\,\, \frac {\partial  F^i}{\partial y} (x, y) = 0\,$ for each
$\,(x, y)
\in E^i$,
\vskip0.5cm
3) $F^i$ has zero values at the endpoints of the arc $E^i$.
\vskip0.5cm

Hence, we can define a function $\tilde{F}$ on the set $\tilde{E}$ as $\tilde{F}(p) = F^i(p)$ for $p \in E^i, i = 1, 2, 3, 4$, and for this function we 
will obviously obtain $\tilde{F} \in C^{2-} (\tilde{E})$ with $ \frac{\partial \tilde{F}}{\partial x} (x, y) = y$ and 
$ \frac{\partial \tilde{F}}{\partial y} (x, y) = 0$ at each point $(x, y) \in \tilde{E}$.
Then, by the classical extension theorem of Whitney (see, for example, Theorem 4 on p. 177 in [S]), there is a function $\overset{\approx}{F} \in C^{2-} 
({\mathbb R}^2_{x, y})$ such that $\overset{\approx}{F} \in C^\infty ({\mathbb R}^2_{x, y} \setminus \tilde{E})$ and 
$\overset{\approx}{F} (p) = \tilde{F} (p)$ for each $p \in \tilde{E}$.
If we restrict the function $\overset{\approx}{F} $ to a disc ${\mathbb D} \subset {\mathbb R}^2_{x, y}$ such that $\tilde{E} \subset {\mathbb D}$ 
and consider a smooth extension of the graph of this restriction to a 2-dimensional sphere $S^2$ embedded into ${\mathbb R}^3_{x, y, z}$, then  the set 
$\overset{\approx}{F} (\tilde{E})$ will be a Jordan curve in $S^2$ of positive 2-dimensional measure and at each point of this curve the tangent plane to 
$S^2$ will coincide with the corresponding plane of the standard contact distribution $\{dz - ydx = 0\}$.

\vskip0.2cm

Now let $G$ be a given strictly pseudoconvex domain in ${\mathbb C}^2$ with $C^\infty$-smooth boundary and let $q$ be a point of $\partial G$.
Then, by the theorem of Darboux, there is a neighbourhood $U$ of $q$ in $\partial G$ and a $C^\infty$-smooth
diffeomorphism $\Phi$ of $U$ onto a neighbourhood $V$ of the origin in ${\mathbb R}^3_{x, y, z}$ such that the distribution of
complex tangencies $\{T^{\mathbb C}_p (\partial G) \}$ will be transformed by $\Phi$ to the standard contact distribution in
${\mathbb R}^3_{x, y, z}$. We can assume without loss of generality that $S^2 \subset V$ (if not, we consider a linear
transformation $x \rightarrow  c x, y \rightarrow c y,  z \rightarrow  c^2 z$ of ${\mathbb R}^3_{x, y, z}$ which preserves the
standard contact structure and use the image of $S^2$ under this transformation with $c > 0$ sufficiently small instead of
$S^2$).  Then $S = \Phi^{-1} (S^2)$ will be a 2-dimensional sphere in $\partial G$ of class $C^{2-}$ and the 
set ${\Gamma} = \Phi^{-1} (\overset{\approx}{F} (\tilde{E})) \subset S$ will be a Jordan curve of positive 2-dimensional measure such that at each 
point 
$p \in {\Gamma}$ the tangent plane $T_pS$ to $S$ is a complex line. 
This proves the second part of the theorem. \bx\\

\bigskip
\parindent 0pt

\end{document}

%% file: fig1.tex
\begin{picture}(0,0)%
\includegraphics{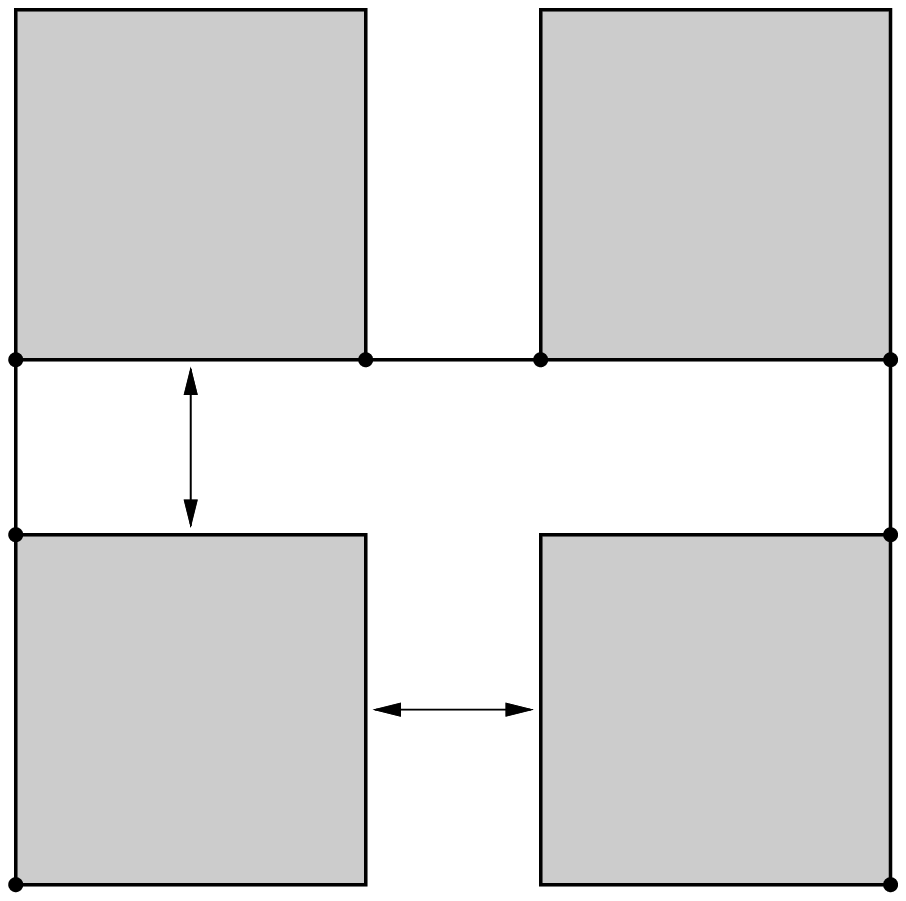}%
\end{picture}%
\setlength{\unitlength}{4144sp}%
\begingroup\makeatletter\ifx\SetFigFont\undefined%
\gdef\SetFigFont#1#2#3#4#5{%
  \reset@font\fontsize{#1}{#2pt}%
  \fontfamily{#3}\fontseries{#4}\fontshape{#5}%
  \selectfont}%
\fi\endgroup%
\begin{picture}(4170,4161)(-84,-3300)
\put(3201,-2431){\makebox(0,0)[b]{\smash{{\SetFigFont{12}{14.4}{\familydefault}{\mddefault}{\updefault}{\color[rgb]{0,0,0}$Q_{3}$}%
}}}}
\put(4071,-1631){\makebox(0,0)[lb]{\smash{{\SetFigFont{12}{14.4}{\familydefault}{\mddefault}{\updefault}{\color[rgb]{0,0,0}$A_{3}$}%
}}}}
\put(4071,-3231){\makebox(0,0)[lb]{\smash{{\SetFigFont{12}{14.4}{\familydefault}{\mddefault}{\updefault}{\color[rgb]{0,0,0}$B_{3}=B_{}=(1,0)$}%
}}}}
\put(3201,-31){\makebox(0,0)[b]{\smash{{\SetFigFont{12}{14.4}{\familydefault}{\mddefault}{\updefault}{\color[rgb]{0,0,0}$Q_{2}$}%
}}}}
\put(2401,-971){\makebox(0,0)[b]{\smash{{\SetFigFont{12}{14.4}{\familydefault}{\mddefault}{\updefault}{\color[rgb]{0,0,0}$A_{2}$}%
}}}}
\put(4071,-831){\makebox(0,0)[lb]{\smash{{\SetFigFont{12}{14.4}{\familydefault}{\mddefault}{\updefault}{\color[rgb]{0,0,0}$B_{2}$}%
}}}}
\put(801,-31){\makebox(0,0)[b]{\smash{{\SetFigFont{12}{14.4}{\familydefault}{\mddefault}{\updefault}{\color[rgb]{0,0,0}$Q_{1}$}%
}}}}
\put(-69,-831){\makebox(0,0)[rb]{\smash{{\SetFigFont{12}{14.4}{\familydefault}{\mddefault}{\updefault}{\color[rgb]{0,0,0}$A_{1}$}%
}}}}
\put(1601,-971){\makebox(0,0)[b]{\smash{{\SetFigFont{12}{14.4}{\familydefault}{\mddefault}{\updefault}{\color[rgb]{0,0,0}$B_{1}$}%
}}}}
\put(801,-2431){\makebox(0,0)[b]{\smash{{\SetFigFont{12}{14.4}{\familydefault}{\mddefault}{\updefault}{\color[rgb]{0,0,0}$Q_{0}$}%
}}}}
\put(-69,-3231){\makebox(0,0)[rb]{\smash{{\SetFigFont{12}{14.4}{\familydefault}{\mddefault}{\updefault}{\color[rgb]{0,0,0}$A_{0}=A_{}=(0,0)$}%
}}}}
\put(-69,-1631){\makebox(0,0)[rb]{\smash{{\SetFigFont{12}{14.4}{\familydefault}{\mddefault}{\updefault}{\color[rgb]{0,0,0}$B_{0}$}%
}}}}
\put(2001,-2291){\makebox(0,0)[b]{\smash{{\SetFigFont{12}{14.4}{\familydefault}{\mddefault}{\updefault}{\color[rgb]{0,0,0}$\alpha$}%
}}}}
\put(871,-1231){\makebox(0,0)[lb]{\smash{{\SetFigFont{12}{14.4}{\familydefault}{\mddefault}{\updefault}{\color[rgb]{0,0,0}$\alpha$}%
}}}}
\end{picture}%

%% file: fig2.tex
\begin{picture}(0,0)%
\includegraphics{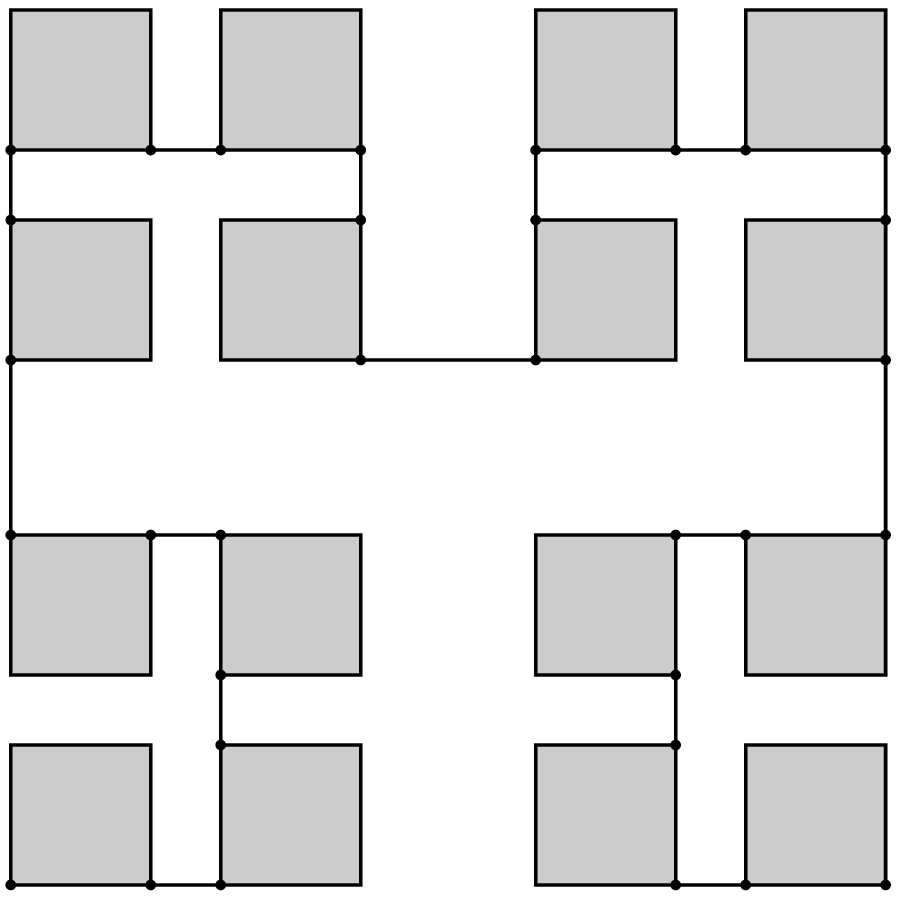}%
\end{picture}%
\setlength{\unitlength}{4144sp}%
\begingroup\makeatletter\ifx\SetFigFont\undefined%
\gdef\SetFigFont#1#2#3#4#5{%
  \reset@font\fontsize{#1}{#2pt}%
  \fontfamily{#3}\fontseries{#4}\fontshape{#5}%
  \selectfont}%
\fi\endgroup%
\begin{picture}(4100,4241)(-49,-3380)
\put(3681,-1951){\makebox(0,0)[b]{\smash{{\SetFigFont{12}{14.4}{\familydefault}{\mddefault}{\updefault}{\color[rgb]{0,0,0}$Q_{30}$}%
}}}}
\put(4036,-1611){\makebox(0,0)[lb]{\smash{{\SetFigFont{12}{14.4}{\familydefault}{\mddefault}{\updefault}{\color[rgb]{0,0,0}{\small$A_{30}=A_{3}$}}%
}}}}
\put(3361,-1511){\makebox(0,0)[b]{\smash{{\SetFigFont{12}{14.4}{\familydefault}{\mddefault}{\updefault}{\color[rgb]{0,0,0}{\small$B_{30}$}}%
}}}}
\put(2721,-1951){\makebox(0,0)[b]{\smash{{\SetFigFont{12}{14.4}{\familydefault}{\mddefault}{\updefault}{\color[rgb]{0,0,0}$Q_{31}$}%
}}}}
\put(3041,-1511){\makebox(0,0)[b]{\smash{{\SetFigFont{12}{14.4}{\familydefault}{\mddefault}{\updefault}{\color[rgb]{0,0,0}{\small$A_{31}$}}%
}}}}
\put(3006,-2151){\makebox(0,0)[rb]{\smash{{\SetFigFont{12}{14.4}{\familydefault}{\mddefault}{\updefault}{\color[rgb]{0,0,0}{\small$B_{31}$}}%
}}}}
\put(2721,-2911){\makebox(0,0)[b]{\smash{{\SetFigFont{12}{14.4}{\familydefault}{\mddefault}{\updefault}{\color[rgb]{0,0,0}$Q_{32}$}%
}}}}
\put(3006,-2671){\makebox(0,0)[rb]{\smash{{\SetFigFont{12}{14.4}{\familydefault}{\mddefault}{\updefault}{\color[rgb]{0,0,0}{\small$A_{32}$}}%
}}}}
\put(3041,-3311){\makebox(0,0)[b]{\smash{{\SetFigFont{12}{14.4}{\familydefault}{\mddefault}{\updefault}{\color[rgb]{0,0,0}{\small$B_{32}$}}%
}}}}
\put(3681,-2911){\makebox(0,0)[b]{\smash{{\SetFigFont{12}{14.4}{\familydefault}{\mddefault}{\updefault}{\color[rgb]{0,0,0}$Q_{33}$}%
}}}}
\put(3361,-3311){\makebox(0,0)[b]{\smash{{\SetFigFont{12}{14.4}{\familydefault}{\mddefault}{\updefault}{\color[rgb]{0,0,0}{\small$A_{33}$}}%
}}}}
\put(4036,-3211){\makebox(0,0)[lb]{\smash{{\SetFigFont{12}{14.4}{\familydefault}{\mddefault}{\updefault}{\color[rgb]{0,0,0}{\small$B_{33}=B_{3}=B_{}$}}%
}}}}
\put(3681,-511){\makebox(0,0)[b]{\smash{{\SetFigFont{12}{14.4}{\familydefault}{\mddefault}{\updefault}{\color[rgb]{0,0,0}$Q_{23}$}%
}}}}
\put(4036,-171){\makebox(0,0)[lb]{\smash{{\SetFigFont{12}{14.4}{\familydefault}{\mddefault}{\updefault}{\color[rgb]{0,0,0}{\small$A_{23}$}}%
}}}}
\put(4036,-811){\makebox(0,0)[lb]{\smash{{\SetFigFont{12}{14.4}{\familydefault}{\mddefault}{\updefault}{\color[rgb]{0,0,0}{\small$B_{23}=B_{2}$}}%
}}}}
\put(3681,449){\makebox(0,0)[b]{\smash{{\SetFigFont{12}{14.4}{\familydefault}{\mddefault}{\updefault}{\color[rgb]{0,0,0}$Q_{22}$}%
}}}}
\put(3361, 49){\makebox(0,0)[b]{\smash{{\SetFigFont{12}{14.4}{\familydefault}{\mddefault}{\updefault}{\color[rgb]{0,0,0}{\small$A_{22}$}}%
}}}}
\put(4036,149){\makebox(0,0)[lb]{\smash{{\SetFigFont{12}{14.4}{\familydefault}{\mddefault}{\updefault}{\color[rgb]{0,0,0}{\small$B_{22}$}}%
}}}}
\put(2721,449){\makebox(0,0)[b]{\smash{{\SetFigFont{12}{14.4}{\familydefault}{\mddefault}{\updefault}{\color[rgb]{0,0,0}$Q_{21}$}%
}}}}
\put(2366,149){\makebox(0,0)[rb]{\smash{{\SetFigFont{12}{14.4}{\familydefault}{\mddefault}{\updefault}{\color[rgb]{0,0,0}{\small$A_{21}$}}%
}}}}
\put(3041, 49){\makebox(0,0)[b]{\smash{{\SetFigFont{12}{14.4}{\familydefault}{\mddefault}{\updefault}{\color[rgb]{0,0,0}{\small$B_{21}$}}%
}}}}
\put(2721,-511){\makebox(0,0)[b]{\smash{{\SetFigFont{12}{14.4}{\familydefault}{\mddefault}{\updefault}{\color[rgb]{0,0,0}$Q_{20}$}%
}}}}
\put(2401,-911){\makebox(0,0)[b]{\smash{{\SetFigFont{12}{14.4}{\familydefault}{\mddefault}{\updefault}{\color[rgb]{0,0,0}{\small$A_{20}=A_{2}$}}%
}}}}
\put(2366,-171){\makebox(0,0)[rb]{\smash{{\SetFigFont{12}{14.4}{\familydefault}{\mddefault}{\updefault}{\color[rgb]{0,0,0}{\small$B_{20}$}}%
}}}}
\put(1281,-511){\makebox(0,0)[b]{\smash{{\SetFigFont{12}{14.4}{\familydefault}{\mddefault}{\updefault}{\color[rgb]{0,0,0}$Q_{13}$}%
}}}}
\put(1636,-171){\makebox(0,0)[lb]{\smash{{\SetFigFont{12}{14.4}{\familydefault}{\mddefault}{\updefault}{\color[rgb]{0,0,0}{\small$A_{13}$}}%
}}}}
\put(1601,-911){\makebox(0,0)[b]{\smash{{\SetFigFont{12}{14.4}{\familydefault}{\mddefault}{\updefault}{\color[rgb]{0,0,0}{\small$B_{13}=B_{1}$}}%
}}}}
\put(1281,449){\makebox(0,0)[b]{\smash{{\SetFigFont{12}{14.4}{\familydefault}{\mddefault}{\updefault}{\color[rgb]{0,0,0}$Q_{12}$}%
}}}}
\put(961, 49){\makebox(0,0)[b]{\smash{{\SetFigFont{12}{14.4}{\familydefault}{\mddefault}{\updefault}{\color[rgb]{0,0,0}{\small$A_{12}$}}%
}}}}
\put(1636,149){\makebox(0,0)[lb]{\smash{{\SetFigFont{12}{14.4}{\familydefault}{\mddefault}{\updefault}{\color[rgb]{0,0,0}{\small$B_{12}$}}%
}}}}
\put(321,449){\makebox(0,0)[b]{\smash{{\SetFigFont{12}{14.4}{\familydefault}{\mddefault}{\updefault}{\color[rgb]{0,0,0}$Q_{11}$}%
}}}}
\put(-34,149){\makebox(0,0)[rb]{\smash{{\SetFigFont{12}{14.4}{\familydefault}{\mddefault}{\updefault}{\color[rgb]{0,0,0}{\small$A_{11}$}}%
}}}}
\put(641, 49){\makebox(0,0)[b]{\smash{{\SetFigFont{12}{14.4}{\familydefault}{\mddefault}{\updefault}{\color[rgb]{0,0,0}{\small$B_{11}$}}%
}}}}
\put(321,-511){\makebox(0,0)[b]{\smash{{\SetFigFont{12}{14.4}{\familydefault}{\mddefault}{\updefault}{\color[rgb]{0,0,0}$Q_{10}$}%
}}}}
\put(-34,-811){\makebox(0,0)[rb]{\smash{{\SetFigFont{12}{14.4}{\familydefault}{\mddefault}{\updefault}{\color[rgb]{0,0,0}{\small$A_{10}=A_{1}$}}%
}}}}
\put(-34,-171){\makebox(0,0)[rb]{\smash{{\SetFigFont{12}{14.4}{\familydefault}{\mddefault}{\updefault}{\color[rgb]{0,0,0}{\small$B_{10}$}}%
}}}}
\put(321,-2911){\makebox(0,0)[b]{\smash{{\SetFigFont{12}{14.4}{\familydefault}{\mddefault}{\updefault}{\color[rgb]{0,0,0}$Q_{00}$}%
}}}}
\put(-34,-3211){\makebox(0,0)[rb]{\smash{{\SetFigFont{12}{14.4}{\familydefault}{\mddefault}{\updefault}{\color[rgb]{0,0,0}{\small$A_{00}=A_{0}=A_{}$}}%
}}}}
\put(641,-3311){\makebox(0,0)[b]{\smash{{\SetFigFont{12}{14.4}{\familydefault}{\mddefault}{\updefault}{\color[rgb]{0,0,0}{\small$B_{00}$}}%
}}}}
\put(1281,-2911){\makebox(0,0)[b]{\smash{{\SetFigFont{12}{14.4}{\familydefault}{\mddefault}{\updefault}{\color[rgb]{0,0,0}$Q_{01}$}%
}}}}
\put(961,-3311){\makebox(0,0)[b]{\smash{{\SetFigFont{12}{14.4}{\familydefault}{\mddefault}{\updefault}{\color[rgb]{0,0,0}{\small$A_{01}$}}%
}}}}
\put(996,-2671){\makebox(0,0)[lb]{\smash{{\SetFigFont{12}{14.4}{\familydefault}{\mddefault}{\updefault}{\color[rgb]{0,0,0}{\small$B_{01}$}}%
}}}}
\put(1281,-1951){\makebox(0,0)[b]{\smash{{\SetFigFont{12}{14.4}{\familydefault}{\mddefault}{\updefault}{\color[rgb]{0,0,0}$Q_{02}$}%
}}}}
\put(996,-2151){\makebox(0,0)[lb]{\smash{{\SetFigFont{12}{14.4}{\familydefault}{\mddefault}{\updefault}{\color[rgb]{0,0,0}{\small$A_{02}$}}%
}}}}
\put(961,-1511){\makebox(0,0)[b]{\smash{{\SetFigFont{12}{14.4}{\familydefault}{\mddefault}{\updefault}{\color[rgb]{0,0,0}{\small$B_{02}$}}%
}}}}
\put(321,-1951){\makebox(0,0)[b]{\smash{{\SetFigFont{12}{14.4}{\familydefault}{\mddefault}{\updefault}{\color[rgb]{0,0,0}$Q_{03}$}%
}}}}
\put(641,-1511){\makebox(0,0)[b]{\smash{{\SetFigFont{12}{14.4}{\familydefault}{\mddefault}{\updefault}{\color[rgb]{0,0,0}{\small$A_{03}$}}%
}}}}
\put(-34,-1611){\makebox(0,0)[rb]{\smash{{\SetFigFont{12}{14.4}{\familydefault}{\mddefault}{\updefault}{\color[rgb]{0,0,0}{\small$B_{03}=B_{0}$}}%
}}}}
\end{picture}%

%% file: fig3.tex
\begin{picture}(0,0)%
\includegraphics{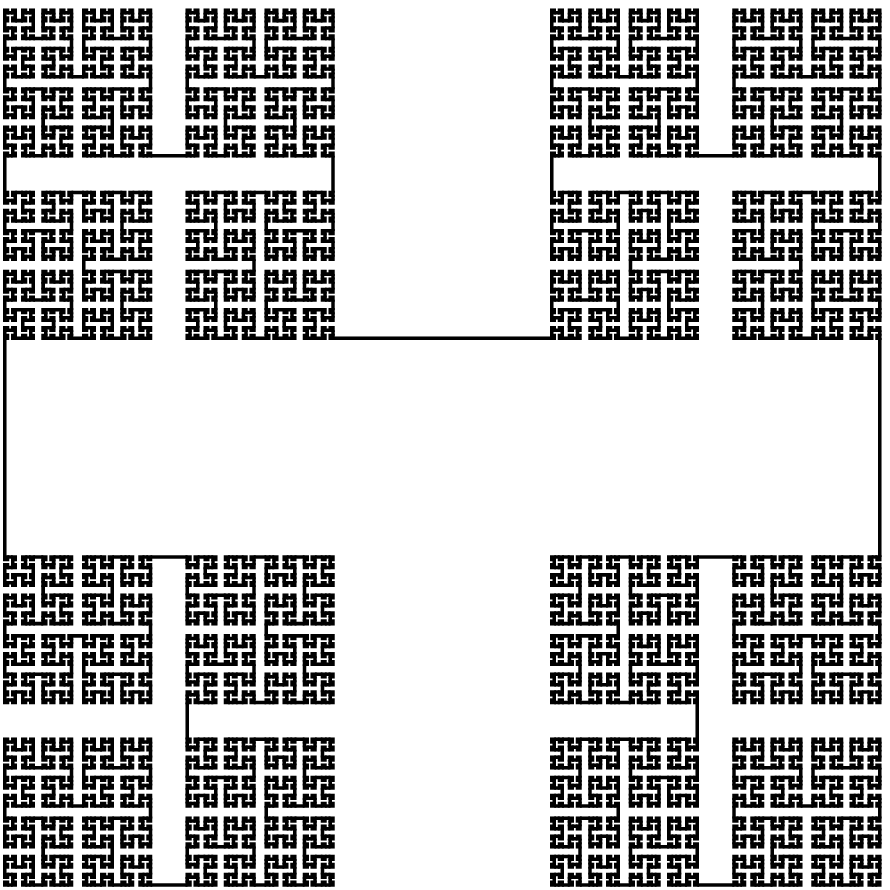}%
\end{picture}%
\setlength{\unitlength}{4144sp}%
\begingroup\makeatletter\ifx\SetFigFont\undefined%
\gdef\SetFigFont#1#2#3#4#5{%
  \reset@font\fontsize{#1}{#2pt}%
  \fontfamily{#3}\fontseries{#4}\fontshape{#5}%
  \selectfont}%
\fi\endgroup%
\begin{picture}(4044,4044)(-21,-3183)
\end{picture}%

%% file: fig5.tex
\begin{picture}(0,0)%
\includegraphics{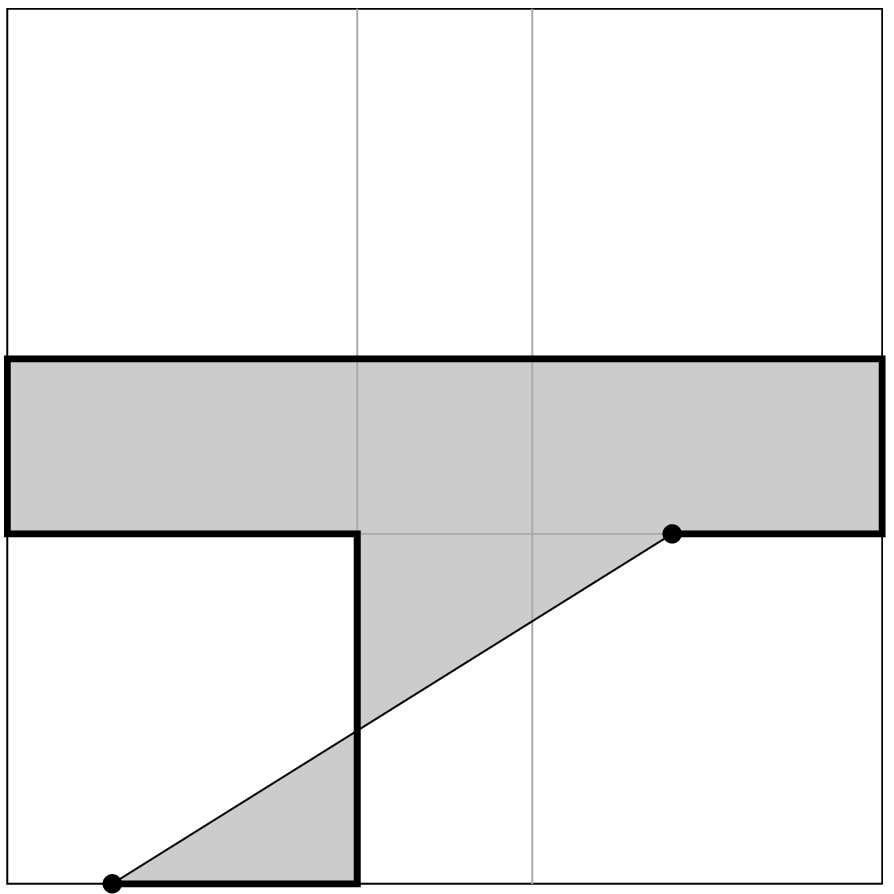}%
\end{picture}%
\setlength{\unitlength}{4144sp}%
\begingroup\makeatletter\ifx\SetFigFontNFSS\undefined%
\gdef\SetFigFontNFSS#1#2#3#4#5{%
  \reset@font\fontsize{#1}{#2pt}%
  \fontfamily{#3}\fontseries{#4}\fontshape{#5}%
  \selectfont}%
\fi\endgroup%
\begin{picture}(4066,4059)(-32,-3208)
\put(1601, 39){\makebox(0,0)[lb]{\smash{{\SetFigFontNFSS{12}{14.4}{\familydefault}{\mddefault}{\updefault}{\color[rgb]{0,0,0}{\Large$Q_{i_1i_2\dots i_{n-1}}$}}%
}}}}
\put(401,-561){\makebox(0,0)[lb]{\smash{{\SetFigFontNFSS{12}{14.4}{\familydefault}{\mddefault}{\updefault}{\color[rgb]{0,0,0}{\Large$J_n^{p,q}$}}%
}}}}
\put(401,-2961){\makebox(0,0)[lb]{\smash{{\SetFigFontNFSS{12}{14.4}{\familydefault}{\mddefault}{\updefault}{\color[rgb]{0,0,0}{\Large$p$}}%
}}}}
\put(3001,-1841){\makebox(0,0)[lb]{\smash{{\SetFigFontNFSS{12}{14.4}{\familydefault}{\mddefault}{\updefault}{\color[rgb]{0,0,0}{\Large$q$}}%
}}}}
\end{picture}%

%% file: fig6b.tex
\begin{picture}(0,0)%
\includegraphics{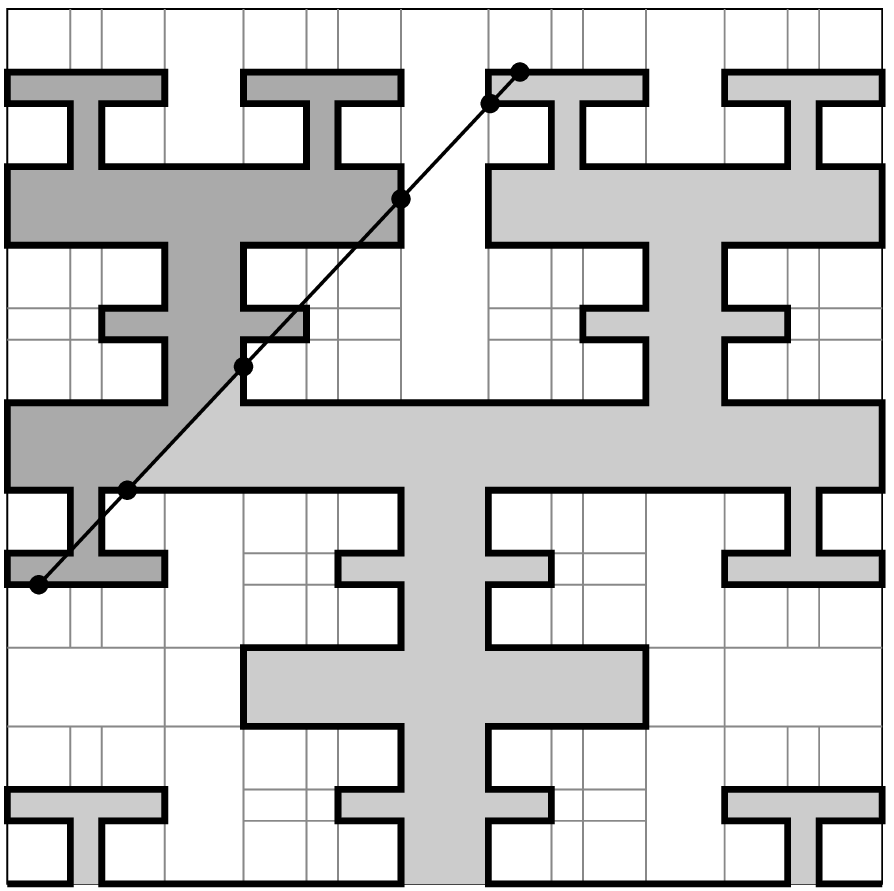}%
\end{picture}%
\setlength{\unitlength}{4144sp}%
\begingroup\makeatletter\ifx\SetFigFontNFSS\undefined%
\gdef\SetFigFontNFSS#1#2#3#4#5{%
  \reset@font\fontsize{#1}{#2pt}%
  \fontfamily{#3}\fontseries{#4}\fontshape{#5}%
  \selectfont}%
\fi\endgroup%
\begin{picture}(4066,4045)(-32,-3194)
\put(1761,-1261){\makebox(0,0)[lb]{\smash{{\SetFigFontNFSS{12}{14.4}{\familydefault}{\mddefault}{\updefault}{\color[rgb]{0,0,0}{\large$\Omega_n^{i_1i_2\dots i_{m_0}}$}}%
}}}}
\put(121,-1261){\makebox(0,0)[lb]{\smash{{\SetFigFontNFSS{12}{14.4}{\familydefault}{\mddefault}{\updefault}{\color[rgb]{0,0,0}{\large$\Sigma_1$}}%
}}}}
\put(121,-2001){\makebox(0,0)[lb]{\smash{{\SetFigFontNFSS{12}{14.4}{\familydefault}{\mddefault}{\updefault}{\color[rgb]{0,0,0}{\large$p=p_0$}}%
}}}}
\put(521,-1561){\makebox(0,0)[lb]{\smash{{\SetFigFontNFSS{12}{14.4}{\familydefault}{\mddefault}{\updefault}{\color[rgb]{0,0,0}{\large$p_0'$}}%
}}}}
\put(1161,-881){\makebox(0,0)[lb]{\smash{{\SetFigFontNFSS{12}{14.4}{\familydefault}{\mddefault}{\updefault}{\color[rgb]{0,0,0}{\large$p_1$}}%
}}}}
\put(1881,-81){\makebox(0,0)[lb]{\smash{{\SetFigFontNFSS{12}{14.4}{\familydefault}{\mddefault}{\updefault}{\color[rgb]{0,0,0}{\large$p_1'$}}%
}}}}
\put(1961,439){\makebox(0,0)[lb]{\smash{{\SetFigFontNFSS{12}{14.4}{\familydefault}{\mddefault}{\updefault}{\color[rgb]{0,0,0}{\large$p_2$}}%
}}}}
\put(2301,639){\makebox(0,0)[lb]{\smash{{\SetFigFontNFSS{12}{14.4}{\familydefault}{\mddefault}{\updefault}{\color[rgb]{0,0,0}{\large$q=p_2'$}}%
}}}}
\end{picture}%

%% file: fig7b.tex
\begin{picture}(0,0)%
\includegraphics{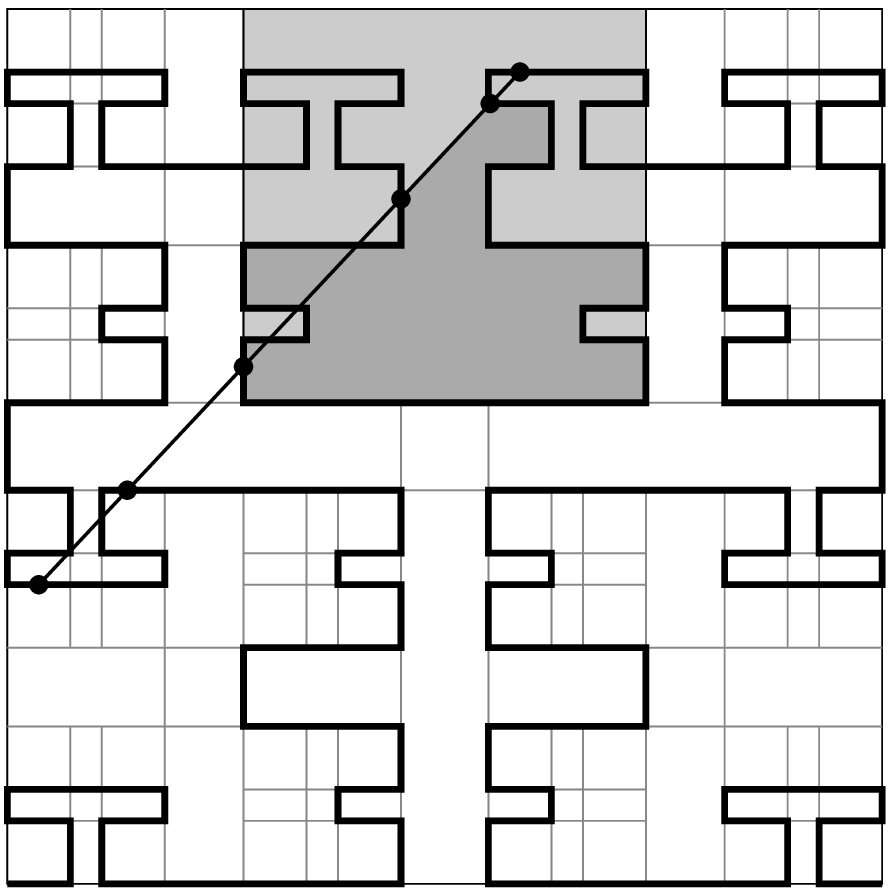}%
\end{picture}%
\setlength{\unitlength}{4144sp}%
\begingroup\makeatletter\ifx\SetFigFontNFSS\undefined%
\gdef\SetFigFontNFSS#1#2#3#4#5{%
  \reset@font\fontsize{#1}{#2pt}%
  \fontfamily{#3}\fontseries{#4}\fontshape{#5}%
  \selectfont}%
\fi\endgroup%
\begin{picture}(4066,4045)(-32,-3194)
\put(1881,-761){\makebox(0,0)[lb]{\smash{{\SetFigFontNFSS{12}{14.4}{\familydefault}{\mddefault}{\updefault}{\color[rgb]{0,0,0}{\large$\Sigma_2$}}%
}}}}
\put(1721,639){\makebox(0,0)[lb]{\smash{{\SetFigFontNFSS{12}{14.4}{\familydefault}{\mddefault}{\updefault}{\color[rgb]{0,0,0}{\large$\Pi$}}%
}}}}
\put(121,-2001){\makebox(0,0)[lb]{\smash{{\SetFigFontNFSS{12}{14.4}{\familydefault}{\mddefault}{\updefault}{\color[rgb]{0,0,0}{\large$p=p_0$}}%
}}}}
\put(521,-1561){\makebox(0,0)[lb]{\smash{{\SetFigFontNFSS{12}{14.4}{\familydefault}{\mddefault}{\updefault}{\color[rgb]{0,0,0}{\large$p_0'$}}%
}}}}
\put(1161,-881){\makebox(0,0)[lb]{\smash{{\SetFigFontNFSS{12}{14.4}{\familydefault}{\mddefault}{\updefault}{\color[rgb]{0,0,0}{\large$p_1$}}%
}}}}
\put(1881,-81){\makebox(0,0)[lb]{\smash{{\SetFigFontNFSS{12}{14.4}{\familydefault}{\mddefault}{\updefault}{\color[rgb]{0,0,0}{\large$p_1'$}}%
}}}}
\put(1961,439){\makebox(0,0)[lb]{\smash{{\SetFigFontNFSS{12}{14.4}{\familydefault}{\mddefault}{\updefault}{\color[rgb]{0,0,0}{\large$p_2$}}%
}}}}
\put(2301,639){\makebox(0,0)[lb]{\smash{{\SetFigFontNFSS{12}{14.4}{\familydefault}{\mddefault}{\updefault}{\color[rgb]{0,0,0}{\large$q=p_2'$}}%
}}}}
\end{picture}%

%% file: fig4a.tex
\begin{picture}(0,0)%
\includegraphics{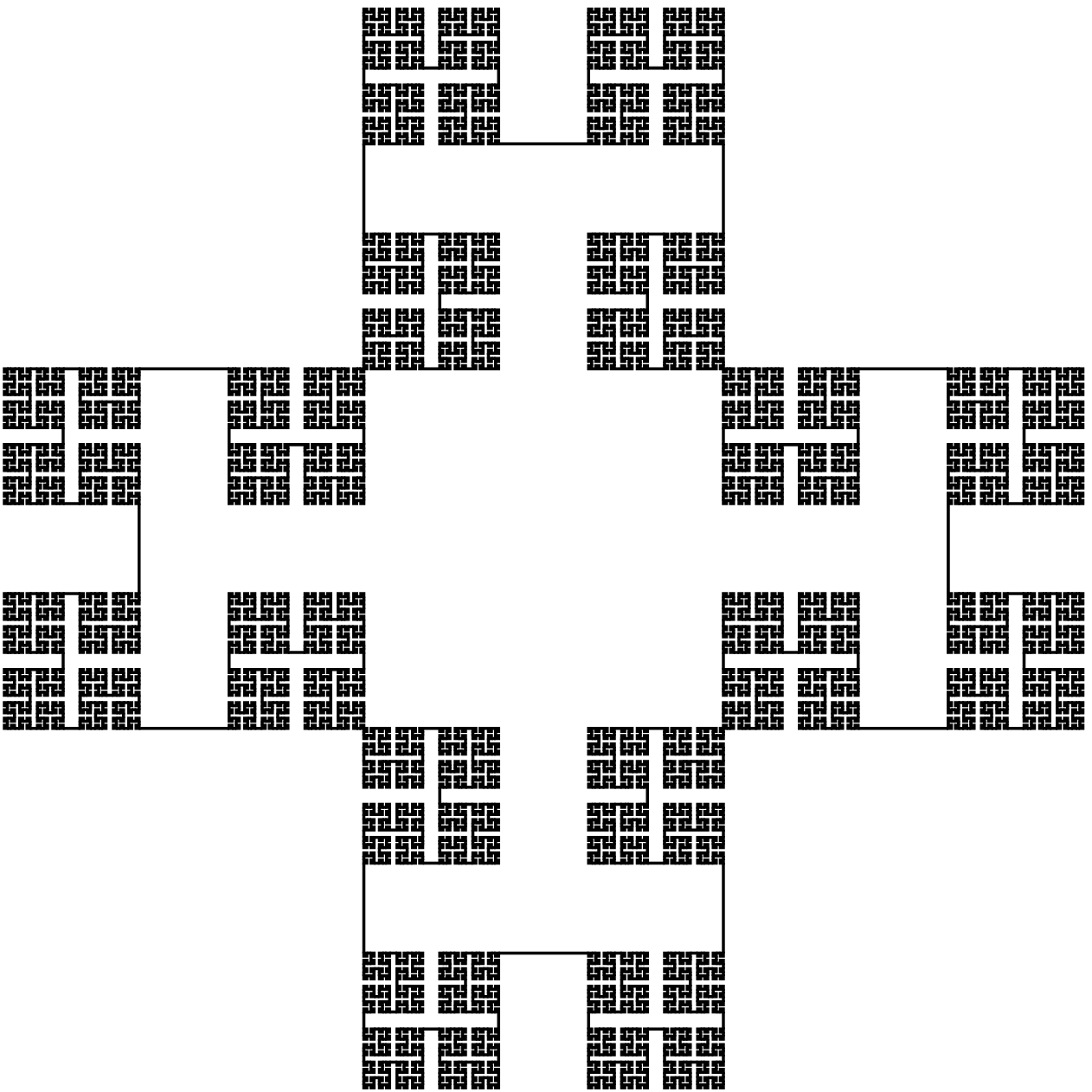}%
\end{picture}%
\setlength{\unitlength}{4144sp}%
\begingroup\makeatletter\ifx\SetFigFontNFSS\undefined%
\gdef\SetFigFontNFSS#1#2#3#4#5{%
  \reset@font\fontsize{#1}{#2pt}%
  \fontfamily{#3}\fontseries{#4}\fontshape{#5}%
  \selectfont}%
\fi\endgroup%
\begin{picture}(5744,5744)(-1921,-4883)
\put(951,-181){\makebox(0,0)[b]{\smash{{\SetFigFontNFSS{12}{14.4}{\familydefault}{\mddefault}{\updefault}{\color[rgb]{0,0,0}$E^1$}%
}}}}
\put(2851,-2081){\makebox(0,0)[b]{\smash{{\SetFigFontNFSS{12}{14.4}{\familydefault}{\mddefault}{\updefault}{\color[rgb]{0,0,0}$E^2$}%
}}}}
\put(-949,-2081){\makebox(0,0)[b]{\smash{{\SetFigFontNFSS{12}{14.4}{\familydefault}{\mddefault}{\updefault}{\color[rgb]{0,0,0}$E^4$}%
}}}}
\put(951,-3981){\makebox(0,0)[b]{\smash{{\SetFigFontNFSS{12}{14.4}{\familydefault}{\mddefault}{\updefault}{\color[rgb]{0,0,0}$E^3$}%
}}}}
\end{picture}%